\documentclass[12,reqno]{amsart}
\usepackage{amsfonts}
\usepackage{subfigure}
\usepackage{graphicx}
\usepackage{amssymb}
\usepackage{amsthm}
\usepackage{amsmath}
\usepackage{dsfont}
\usepackage{a4wide}
\usepackage{pgf,pgfarrows,pgfautomata,pgfheaps,pgfnodes,pgfshade}
\usepackage{url}

\usepackage[numbers]{natbib}

\usepackage{hyperref}
\numberwithin{equation}{section}

\newtheorem{defi}{defi}[section]
\newtheorem{thm}[defi]{Theorem}
\newtheorem{lemm}[defi]{Lemma}
\newtheorem{remark}[defi]{Remark}

\newtheorem{assum}[defi]{Assumption}
\newtheorem{prop}[defi]{Proposition}

\begin{document}

	\title[ ALM with Epstein-Zin recursive utility under partial information]{Asset-liability management with Epstein-Zin utility under stochastic interest rate and unknown market price of risk}
	\author[W. Kuissi-Kamdem]{Wilfried Kuissi-Kamdem}
	\address{African Institute for Mathematical Sciences, Ghana}
	\address{Department of Mathematics, University of Rwanda, Rwanda}
	\address{Department of Mathematical Stochastics, University of Freiburg, Germany}
	\email{donatien@aims.edu.gh, wilfried.kuissi.kamdem@stochastik.uni-freiburg.de}
	\thanks{This work was supported by a grant from the African Institute for Mathematical Sciences, with financial support from the Government of Canada, provided through Global Affairs Canada, and the International Development Research Centre. I also gratefully acknowledge support from the Carl Zeiss Foundation in the Rescale project.}
	
	\date{}

	
	\maketitle
	
	\begin{abstract}
		We study a continuous-time consumption--investment problem with Epstein--Zin recursive utility under partial information, where the market price of risk is unobservable and the investor faces a terminal liability. The introduction of terminal liabilities fundamentally changes the recursive utility optimisation problem by requiring a new ansatz for the value process. This reformulation leads to a coupled linear forward--backward stochastic differential equation (FBSDE) with unbounded random coefficients, which lies outside the scope of standard well-posedness results for coupled FBSDEs. Exploiting the special structure of the equation, we develop a decoupling reduction method that yields an explicit solution under mild assumptions on the model coefficients. We further establish the Malliavin differentiability of the solution, allowing the optimal investment strategy to be represented in terms of conditional expectations involving Malliavin derivatives, which can be efficiently approximated using Monte Carlo methods. The resulting framework provides explicit expressions for the optimal consumption strategy, portfolio allocation, and value function in a stochastic volatility market with stochastic interest rates and an unobservable market price of risk, extending the classical Kim--Omberg model. Finally, we quantify the utility loss arising from ignoring learning about the market price of risk, thereby illustrating the economic significance of partial information in asset--liability management.
	\end{abstract}


\section{Introduction}
Asset--liability management (ALM) plays a central role in the investment decisions of banks, insurance companies, pension funds, and individual investors; see \cite{luckner2003asset}. In many applications, investors seek consumption and investment strategies that account not only for the evolution of asset values but also for future financial obligations. While continuous-time ALM problems have been extensively studied under mean--variance criteria and time-additive utility preferences, comparatively little is known when preferences are described by Epstein--Zin recursive utility, particularly under partial information.

Epstein--Zin utility has become a standard framework in intertemporal portfolio choice because it disentangles the coefficient of relative risk aversion from the elasticity of intertemporal substitution, thereby overcoming an important limitation of classical time-additive utility. Since its introduction by \cite{epstein1989substitution}, recursive utility has been widely used in portfolio optimisation and asset pricing. Existing studies have developed analytical tools for solving consumption--investment problems under recursive preferences, including stochastic filtering, forward--backward stochastic differential equations, and Malliavin calculus.

There is now enough evidence in the literature that stock returns are predictable; see \cite{ait2001variable} for a survey. In \cite{xia2001learning}, the predictive variables were assumed to be unobservable, and since then this framework has been widely adopted in the portfolio optimisation literature. However, in contrast to the rich literature on partial-information portfolio optimisation under classical time-additive utility preferences (see, e.g., \cite{ceci2025portfolio,escobar2016portfolio,frey2012portfolio} for a review), relatively little attention has been devoted to recursive utility maximisation under partial information. Notable exceptions are \cite{chen2021dynamic,ma2023consumption}, which investigate Epstein--Zin utility maximisation under partial information in different infinite-horizon settings. To the best of our knowledge, neither study incorporates terminal liabilities into the optimisation problem.

The present paper studies an Epstein--Zin consumption--investment problem under partial information in which the investor is exposed to a terminal liability and the market price of risk is unobservable. Although each of these ingredients has been studied separately, their combination gives rise to a fundamentally different optimisation problem. In particular, the presence of the terminal liability necessitates a new ansatz for the optimal value process, leading to a coupled linear FBSDE with unbounded random coefficients. To the best of our knowledge, the resulting FBSDE falls outside the scope of existing well-posedness theories, which typically rely on monotonicity and/or Lipschitz conditions that are not satisfied in our setting. Consequently, the mathematical challenge is not merely the inclusion of a liability term but the development of an analytical framework capable of solving the resulting optimisation problem.

From an application perspective, our work provides the first analysis of an Epstein--Zin asset--liability management problem under partial information. The market model extends the classical Kim--Omberg framework \cite{kim1996dynamic} by simultaneously incorporating stochastic interest rates, stochastic volatility, and an unobservable stochastic market price of risk. The resulting explicit characterisation of the optimal consumption and investment strategies allows us to investigate the economic value of learning and to quantify the welfare loss incurred when investors ignore information about the latent market price of risk.

The main contributions of this paper are as follows.
\begin{enumerate}	
	\item We formulate, to the best of our knowledge, the first continuous-time consumption--investment problem with Epstein--Zin recursive utility, terminal liabilities, and partial information. Our framework allows the terminal liability to depend on the entire path of the risky assets, thereby encompassing a broad class of path-dependent contingent claims, including equity-linked securities and convertible bonds. While recent studies \cite{dadzie2025robust,kuissi2026optimal} have considered liabilities under recursive preferences, they are restricted to full-information settings with substantially simpler market dynamics.
	
	\item A key technical contribution is that the presence of terminal liabilities fundamentally changes the recursive utility optimisation problem. In particular, it necessitates a new ansatz for the optimal value process, which differs from those employed in the existing literature on Epstein--Zin utility maximisation. This representation leads to a coupled linear FBSDE with unbounded random coefficients. To the best of our knowledge, this class of FBSDEs has not previously arisen in recursive utility optimisation.
	
	\item The coupled FBSDE generated by our optimisation problem falls outside the scope of existing well-posedness results, including those in \cite{ankirchner2022transformation,hu1995solution,ma2015well,pardoux1990adapted,pardoux1999forward,peng1999fully,xie2020exploration}, which typically require monotonicity and/or Lipschitz conditions that are not satisfied in our setting. To overcome this difficulty, we develop a decoupling reduction method that exploits the special structure of the equation. In contrast to the approaches of \cite{ma1994solving,yong1999linear,yong2006linear}, where coupled linear FBSDEs are reduced to backward stochastic differential equations, our construction reduces the problem to a system of forward stochastic differential equations. This reduction is structurally different and establishes explicit solvability.
	
	\item The proposed methodology yields explicit expressions for the optimal consumption strategy, portfolio allocation, and value function. Under mild additional assumptions on the liability, we establish the Malliavin differentiability of the solution, allowing the optimal investment strategy to be represented through conditional expectations involving Malliavin derivatives. This representation is particularly attractive from a computational perspective, as it can be efficiently approximated using Monte Carlo simulation.
	
	\item In the special case of the classical CRRA utility, where the coefficient of relative risk aversion and the elasticity of intertemporal substitution are inversely related, the optimisation problem \eqref{Main problem} is related to \cite{horst2014forward}. Our framework, however, incorporates intermediate consumption and admits a simpler solution methodology. Specifically, \cite{horst2014forward} characterises the optimal investment strategy through a coupled nonlinear FBSDE, for which no explicit solution is available. By contrast, we characterise the optimal consumption--investment pair by a coupled linear FBSDE and derive an explicit solution. This explicit characterisation not only facilitates the theoretical analysis but also makes our results readily amenable to numerical implementation.
	\item We quantify the utility loss incurred by investors who ignore the opportunity to learn about the market price of risk; see Section \ref{Welfare loss}. Following \cite{escobar2016portfolio}, we measure this loss in terms of the percentage of initial wealth that investors would be willing to forgo in exchange for access to learning, commonly referred to as the \emph{welfare loss}. Our numerical results show that, in the presence of liabilities, the welfare loss is increasing in the investors' initial wealth, whereas it is independent of initial wealth when liabilities are absent. Furthermore, the welfare loss decreases with the coefficient of relative risk aversion and increases with the elasticity of intertemporal substitution. We also find that the welfare loss is smaller under CRRA utility than under Epstein--Zin preferences.
\end{enumerate}

The remainder of the present paper is structured as follows. We introduce the model and formulate the problem in Section \ref{Setting}. In Section \ref{Main results} we give the main results of this paper. Finally, in Section \ref{Welfare loss} we determine the utility loss and perform some numerical analysis. Some additional technical proofs are postponed in the Appendix.

\section{Model and problem formulation}\label{Setting}
We consider a probability space $(\Omega,\mathbb{F},\mathbb{P})$ equipped with the complete filtration $(\mathcal{F}_{t})_{t\in[0,T]}$ supporting a standard three dimensional Wiener process $W_:=(W^{1},W^{2},W^{3})$.

\subsection{The financial market}\label{The financial market}
The market dynamics are described by an extension of the Kim--Omberg model \cite{kim1996dynamic} that accommodates stochastic interest rates and an unobservable stochastic market price of risk. Specifically, we consider a dynamic financial environment with three traded assets and one non-traded financial index. The traded assets consist of one money market account $S^{0}$, one stock $S$ and one zero-coupon bond $B$ maturing at time $T$. The money market account follows
\begin{align}\label{Money market account}
	\mathrm{d}S_{t}^{0}&=r_{t}S_{t}^{0}\mathrm{d}t,~S_{0}^{0}>0,
\end{align}
where the stochastic short-term interest rate $r$ follows the Vasicek model \cite{vasicek1977equilibrium}; that is, $r$ is an Ornstein-Uhlenbeck (OU) process given by
\begin{align}\label{Stochastic interest rate}
	\mathrm{d}r_{t}&=\kappa_{r}\big(\mu_{r}-r_{t}\big)\mathrm{d}t+\sigma_{r}\Big(\rho_{rS}\mathrm{d}W_{t}^{1}+\sqrt{1-\rho_{rS}^{2}}\mathrm{d}W_{t}^{2}\Big),
\end{align}
with correlation coefficient $\rho_{rS}\in(-1,1)$, speed of mean reversion $\kappa_{r}$, long run mean $\mu_{r}$ and volatility $\sigma_{r}>0$.

Following \cite{vasicek1977equilibrium}, the zero-coupon bond evolves according to
\begin{align}\label{Zero-coupon bond}
	\mathrm{d}B_{t}=B_{t}\Big(\big(r_{t}+\mu_{B}(t)\big)\mathrm{d}t+\sigma_{B}(t)\Big(\rho_{rS}\mathrm{d}W_{t}^{1}+\sqrt{1-\rho_{rS}^{2}}\mathrm{d}W_{t}^{2}\Big)\Big),
\end{align}
with correlation coefficient $\rho_{rS}$, excess rerun of the bond $\mu_{B}(t):=\phi_{B}\sigma_{B}(t)$ and volatility $\sigma_{B}(t):=\sigma_{r}\frac{1-\exp(-\kappa_{r}(T-t))}{\kappa_{r}}$. We assume that the investors follow a roll-over strategy for the bond investment and keeps the maturity of the bond in their portfolio constant. This is a common assumption in the literature on portfolio choice with stochastic interest rates; see \cite{escobar2016portfolio} and reference therein.

The stock price has dynamics given by
\begin{align}\label{Stock}
	\mathrm{d}S_{t}=S_{t}\Big(\big(r_{t}+\beta\sigma_{t}R_{t}\big)\mathrm{d}t+\sigma_{t}\mathrm{d}W_{t}^{1}\Big),~S_{0}>0,
\end{align}
with $\sigma$ an uniformly positive process adapted to the filtration generated by $r$ and $S$, and $\beta\ne0$. Similar setup is also considered in \cite{detemple2003monte} under full information where the expected stock return depends only on the interest rate.

In \eqref{Stock}, $R$ is an $\mathbb{R}$-valued non-traded financial index which follows a linear mean-reverting dynamics given by
\begin{align}\label{Risk premium}
	\mathrm{d}R_{t}&=\kappa_{R}\big(\mu_{R}-R_{t}\big)\mathrm{d}t+\sigma_{R}\Big(\rho_{RS}\mathrm{d}W_{t}^{1}+\rho_{Rr}\mathrm{d}W_{t}^{2}+\sqrt{1-\rho_{RS}^{2}-\rho_{Rr}^{2}}\mathrm{d}W_{t}^{3}\Big),
\end{align}
with correlation coefficients $\rho_{RS},\rho_{Rr}\in[-1,1]$, speed of mean reversion $\kappa_{R}$, long run mean $\mu_{R}$ and volatility $\sigma_{R}>0$. In the sequel, following \cite{escobar2016portfolio}, we assume that $\rho_{Rr}:=\frac{\alpha-\rho_{rS}\rho_{RS}}{\sqrt{1-\rho_{rS}^{2}}}$ for $\alpha\in\mathbb{R}$ such that $\rho_{RS}^{2}+\rho_{Rr}^{2}\in[0,1]$. Hence, the process $(R_{t})_{t\in[0,T]}$ plays the role of the market price of risk.

Hence, investors choose the consumption rate $c_{t},~t\in[0,T]$, and the amounts $\pi_{t}^{S}$ and $\pi_{t}^{B}$ to be invested in the stock and in the bond, respectively. For such $(c,\pi^{S},\pi^{B})$, the wealth process $X$ of the investors with initial endowment $x$ at time $0$ evolves according to the stochastic differential equation
\begin{align}\label{Wealth}
	\mathrm{d}X_{t}&=\Big(r_{t}X_{t}+\pi_{t}^{S}\beta\sigma_{t}R_{t}+\pi_{t}^{B}\phi_{B}\sigma_{B}(t)\Big)\mathrm{d}t+\Big(\pi_{t}^{S}\sigma_{t}+\pi_{t}^{B}\sigma_{B}(t)\rho_{rS}\Big)\mathrm{d}W_{t}^{1}\notag\\
	&\phantom{X}+\pi_{t}^{B}\sigma_{B}(t)\sqrt{1-\rho_{rS}^{2}}\mathrm{d}W_{t}^{2}-c_{t}\mathrm{d}t.
\end{align}
Note that the market is incomplete (the number of traded assets being less than the number of Wiener processes).

\subsection{The partial information framework}\label{The partial information framework}
We assume that the risk premium $R_{t},~t\in[0,T]$, is not directly observable by the investors. Hence, the investors have no direct information on the return of the stock. The available information flow comes from past realisations/observation of two processes: the stochastic interest rate $r$ and the stock $S$. We
introduce the observation filtration as $\mathbb{F}^{r,S}:=\mathbb{F}^{r}\vee\mathbb{F}^{S}$, with $\mathbb{F}^{r}:=(\mathcal{F}_{t}^{r})_{t\in[0,T]}$ and $\mathbb{F}^{S}:=(\mathcal{F}_{t}^{S})_{t\in[0,T]}$ being the natural filtration of $r$ and $S$, respectively. We assume that $\mathbb{F}^{r,S}$ is completed with $\mathbb{P}$-null sets and right-continuous.

We end this section with the definition of some spaces that are used throughout. Let $\mathcal{C}$ be the set of $\mathbb{F}^{r,S}$-non-negative progressively measurable processes on $[0,T]\times\Omega$. For $c\in\mathcal{C}$ and $t<T$, $c_{t}$ denotes the consumption rate at time $t$ and $c_{T}$ represents a lumpsum consumption at the finite time horizon $T$. Let $\mathcal{L}_{\mathbb{P}}^{q},~{q>1}$, denotes the space of $\mathcal{F}_{T}^{r,S}$-measurable $\mathbb{R}$-valued random variables $X$ such that $\mathbb{E}[|X|^{q}~\big|\mathcal{F}_{0}^{r,S}]<\infty$. Let $\mathcal{H}_{\mathbb{P}}^{q},~{q>1}$, denotes the space of $\mathbb{F}^{r,S}$-predictable $\mathbb{R}$-valued processes $(Y_{t})_{t\in[0,T]}$ such that $\mathbb{E}[\int_{0}^{T}|Y_{t}|^{q}\mathrm{d}t~\big|\mathcal{F}_{0}^{r,S}]<\infty$. {Let $\mathcal{S}_{\mathbb{P}}^{q},~q>1$, denotes the space of continuous $\mathbb{F}^{r,S}$-adapted $\mathbb{R}$-valued processes $(Y_{t})_{t\in[0,T]}$ such that $\mathbb{E}[\sup_{t\le T}|Y_{t}|^{q}~\big|\mathcal{F}_{0}^{r,S}]<\infty$.} Let $\mathbb{H}_{\mathbb{P}}^{q},~{q>1}$, denotes the space of $\mathbb{F}^{r,S}$-predictable $\mathbb{R}^{2}$-valued processes $(Z_{t})_{t\in[0,T]}$ such that $\mathbb{E}[(\int_{0}^{T}|Z_{t}|^{2}\mathrm{d}t)^{\frac{q}{2}}~\big|\mathcal{F}_{0}^{r,S}]$ $<\infty$. Note that similar spaces can and will be defined under another probability measure $\mathbb{Q}$, by replacing $\mathbb{P}$ with $\mathbb{Q}$ in the subscripts of the corresponding spaces, and taking expectations with respect to $\mathbb{Q}$.

\subsection{The Epstein-Zin utility maximisation problem with partial information}\label{The recursive utility maximisation problem}
An agent's preference over $\mathcal{C}$-valued consumption is given by the Epstein-Zin recursive preference. To describe this preference, let $\delta>0$ represent the discounting rate, $0<\gamma\neq1$ be the relative risk aversion, and $0<\psi\neq1$ be the elasticity of intertemporal substitution coefficient (EIS). Then, the Epstein–Zin aggregator is defined by
\begin{align}\label{Epstein-Zin generator}
	f(c,v)&:=\delta e^{-\delta t} \frac{c^{1-\frac{1}{\psi}}}{1-\frac{1}{\psi}}((1-\gamma)v)^{1-\frac{1}{\theta}},\text{ with }~\theta:=\frac{1-\gamma}{1-\frac{1}{\psi}},
\end{align}
and the bequest utility function by $h(c):=e^{-\delta\theta T}\frac{c^{1-\gamma}}{1-\gamma}$. Hence, the Epstein-Zin utility over the consumption stream $c\in\mathcal{C}$ on a finite time horizon $T$ is a process $V^{c}$ which satisfies
\begin{align}\label{Epstein-Zin utility}
	V_{t}^{c}&=\mathbb{E}\Big[h(c_{T})+\int_{t}^{T}f(c_{s},V_{s}^{c})\mathrm{d}s~\big|\mathcal{F}_{t}^{r,S}\Big]~\text{ for }t\in[0,T].
\end{align}

We focus on the empirically relevant case in which both the coefficient of relative risk aversion and the elasticity of intertemporal substitution are greater than one; i.e., $\gamma>1,\psi>1$. Indeed, there are at least two strands of literature that provide support in the context of consumption and portfolio choice. First, in order to address asset pricing anomalies, such as the excess volatility puzzle and the credit risk puzzle, the authors in \cite{benzoni2011explaining,bhamra2010levered} require $\psi>1$ to align their models with observed data. The empirical findings also indicate that risk aversion satisfies $\gamma>1$; see \cite[on p.228]{xing2017consumption} for a further detailed discussion of the parameter choices for $\gamma$ and $\psi$. More general parameter configurations are analysed by \cite{matoussi2018convex} in a finite-horizon setting, and by \cite{herdegen2023infiniteI,shigeta2026economic} in an infinite-horizon framework, both in a full-information framework.
\begin{remark}\label{Special CRRA utility}
	While our analysis focuses on the empirically relevant parameter regime $\gamma>1,\,\psi>1$, the analysis developed in this paper also applies under the parameter configuration $\gamma>1$ and $\gamma\psi=1$, the latter corresponding to the special case of CRRA utility. The only point at which the CRRA case requires separate treatment is in the proof of Theorem \ref{Mainresult1}, where we explicitly indicate the necessary modifications.
\end{remark}
\begin{defi}\label{Admissible consumption}
	A consumption stream $c\in\mathcal{C}$ is said to be \textit{admissible} if Equation \eqref{Epstein-Zin utility} admits a unique solution $V^{c}$ within the class of processes of class (D) satisfying $(1-\gamma)V^{c}>0$. The set of all admissible consumption streams is denoted by $\mathcal{C}_{a}$.
\end{defi}
The set $\mathcal{C}_{a}$ in Definition~\ref{Admissible consumption} is similar to the one considered in \cite{herdegen2023infiniteI,matoussi2018convex}. All known sufficient conditions for the existence of Epstein–Zin utility over a finite time horizon are summarised in \cite[Prop.~2.1]{matoussi2018convex}, which, in particular, ensures that $\mathcal{C}_{a}\ne\emptyset$.

In the present paper, we investigate a consumption–investment problem with Epstein--Zin utility under terminal liability constraint and partial information. Specifically, we consider an investor with a terminal liability $K_{T}$ at time $T$ and Epstein--Zin recursive preferences, and assume that the market price of risk is unobservable. The liability $K_{T}$ is allowed to take either positive or negative values. The objective is to determine an optimal consumption and investment strategy $(c^{*},\pi^{S,*},\pi^{B,*})$ that solves the following optimisation problem:
\begin{align}\label{Main problem}
	\widehat{V}_{0}:=\sup_{(c,\pi^{S},\pi^{B})\in\mathcal{P}}~\mathbb{E}\Big[h(X_{T}-K_{T})+\int_{0}^{T}f(c_{t},V_{t}^{c})\mathrm{d}t~\big|\mathcal{F}_{0}^{r,S}\Big],
\end{align}
where $\mathcal{P}$ is a subset of the set of $\mathbb{R}^{3}$-valued $\mathbb{F}^{r,S}$-adapted processes. A precise definition of the set $\mathcal{P}$ is postponed in Definition \ref{Permissible strategies}.

A key feature of the stochastic optimisation problem \eqref{Main problem} is that the supremum is taken over strategies adapted to the observation filtration $\mathbb{F}^{r,S}$, rather than the global filtration $\mathbb{F}$. This places us in the setting of stochastic optimisation under partial information. To address this challenge, we follow the approach of \cite{fleming1982optimal} and introduce an auxiliary separated problem. In the separated formulation, all state variables are adapted to $\mathbb{F}^{r,S}$. Establishing this requires tools from stochastic filtering theory, which will be presented in Section \ref{Reduction to the observable filtration}. See \cite{liptser2013statisticsI} for more details on the subject.

\section{Main results}\label{Main results}
\subsection{Reduction to the observable filtration}\label{Reduction to the observable filtration}
Mathematically the financial market is described in terms of a partially observable triple of processes $(R,r,S)$, where $R$ is called the unobservable signal, and $r$ and $S$ the observation processes.To reduce the problem to the observation filtration, we consider the process $(m_{t})_{t\in[0,T]}$ and the function $v:[0,T]\to\mathbb{R}_{+}$ defined by
\begin{align}\label{Conditional characteristics of R}
	m_{t}&:=\mathbb{E}\big[R_{t}~|\mathcal{F}_{t}^{r,S}\big]~\text{ and }~v_{t}:=\mathbb{E}\big[\big(R_{t}-m_{t}\big)^{2}~|\mathcal{F}_{t}^{r,S}\big]~\text{ for }~t\in[0,T].
\end{align}
Then, we obtain the following results.
\begin{prop}\label{Conditional mean-variance pair}
	Let the conditional mean-variance pair $(m_{t},v_{t})_{t\in[0,T]}$ be defined as in \eqref{Conditional characteristics of R}. Then, $(m_{t},v_{t}),~t\in[0,T]$, solves the system
	\begin{align*}
		\begin{cases}
			\mathrm{d}m_{t}&=\kappa_{R}\big(\mu_{R}-m_{t}\big)\mathrm{d}t+\big(\sigma_{R}\rho_{RS}+\beta v_{t}\big)\mathrm{d}I_{t}^{1}+\big(\sigma_{R}\rho_{Rr} -\rho_{rS}\beta(1-\rho_{rS}^{2})^{-\frac{1}{2}}v_{t}\big)\mathrm{d}I_{t}^{2}\\
			\mathrm{d}v_{t}&=\Big(\sigma_{R}^{2}-2\kappa_{R}v_{t}-\big(\sigma_{R}\rho_{RS}+\beta v_{t}\big)^{2}-\big(\sigma_{R}\rho_{Rr} -\rho_{rS}\beta(1-\rho_{rS}^{2})^{-\frac{1}{2}}v_{t}\big)^{2}\Big)\mathrm{d}t,
		\end{cases}
	\end{align*}
	with $m_{0}=\mathbb{E}\big[R_{0}~|\mathcal{F}_{t}^{r,S}\big],v_{0}=\mathbb{E}\big[\big(R_{0}-m_{0}\big)^{2}~|\mathcal{F}_{t}^{r,S}\big]$ and $I=(I_{t}^{1},I_{t}^{2})_{t\in[0,T]}$ the $\mathbb{R}^{2}$-valued process, called the \textit{innovation process}, given by
	\begin{align}\label{Innovation process}
		I_{t}^{1}&:=W_{t}^{1}+\beta\int_{0}^{t}\big(R_{s}-m_{s}\big)\mathrm{d}s,~~I_{t}^{2}:=W_{t}^{2}-\frac{\rho_{rS}}{(1-\rho_{rS}^{2})^{1/2}}\beta\int_{0}^{t}\big(R_{s}-m_{s}\big)\mathrm{d}s
	\end{align}
	is a two dimensional Brownian motion under the filtration $\mathbb{F}^{r,S}$ and the probability $\mathbb{P}$.
\end{prop}
\begin{proof}
	The proof follows similar arguments as in the proof of proposition $1$ in \cite{escobar2016portfolio} for $\sigma_{\lambda},\kappa_{\lambda},\rho_{S\lambda},\hat{\rho}_{\lambda},\hat{\rho}_{\lambda P}$ and $\hat{\rho}_{\lambda\beta}$ therein substituted by $0,0,0,1,0$ and $0$, respectively.
\end{proof}
Using the definition of the innovation process, given by \eqref{Innovation process}, we can equivalently write the dynamics of the interest rate $(r_{t})_{t\in[0,T]}$ and the wealth process $(X_{t})_{t\in[0,T]}$ as follows:
\begin{align}
	\mathrm{d}r_{t}&=\kappa_{r}\big(\mu_{r}-r_{t}\big)\mathrm{d}t+\sigma_{r}\Big(\rho_{rS}\mathrm{d}I_{t}^{1}+\sqrt{1-\rho_{rS}^{2}}\mathrm{d}I_{t}^{2}\Big)\label{Filtered interest rate}\\
	\mathrm{d}X_{t}&=\Big(r_{t}X_{t}+\pi_{t}^{\top}\eta_{t}\Big)\mathrm{d}t+\pi_{t}^{\top}\mathrm{d}I_{t}-c_{t}\mathrm{d}t,~X_{0}=x,\label{Filtered wealth}
\end{align}
where $\Sigma_{t}:=\begin{pmatrix}
	\sigma_{t} & 0\\ \sigma_{B}(t)\rho_{rS} & \sigma_{B}(t)\sqrt{1-\rho_{rS}^{2}}
\end{pmatrix}$, $\mu_{t}:=\begin{pmatrix}
	\beta\sigma_{t}m_{t}\\ \phi_{B}\sigma_{B}(t)
\end{pmatrix}$, $\pi_{t}^{\top}:=(\pi_{t}^{S},\pi_{t}^{B})\Sigma_{t}$ and $\eta_{t}:=\Sigma_{t}^{-1}\mu_{t}=\big(\beta m_{t},~\big(1-\rho_{rS}^{2}\big)^{-\frac{1}{2}}\big(-\beta\rho_{rS}m_{t}+\phi_{B}\big)\big)^{\top}$ for $t\in[0,T]$. In the sequel, $X^{c,\pi}$ denotes the solution of \eqref{Filtered wealth}.
\begin{defi}[Admissible strategies]
	A pair $(c,\pi)$ of $\mathbb{F}^{r,S}$-adapted consumption-portfolio strategy is \textit{admissible} if \eqref{Filtered wealth} admits a unique strong solution. We denote by $\mathcal{A}$ the set of admissible consumption-portfolio strategies.
\end{defi}
Note that in \eqref{Filtered wealth} the unobservable market price of risk process $(R_{t})_{t\in[0,T]}$ does not appear anymore, and all coefficients are adapted to the observation filtration $\mathbb{F}^{r,S}$. From now on, all (conditional) expectations are taken under the observation filtration and all SDEs are under the Brownian motion $(I_{t})_{t\in[0,T]}$. Moreover, we still denote by $K_{T}$ the filtered value of the liability with respect to $\mathcal{F}_{T}^{r,S}$.

\subsection{Solution to the optimisation problem}\label{Solution to the optimisation problem}
The presence of the terminal liability $K_T$ prevents a direct application of the standard homothetic ansatz or \textit{separation of variables} property commonly used in the Epstein--Zin portfolio optimisation literature, see, e.g., \cite{aurand2023epstein,feng2026consumption,kraft2017optimal,xing2017consumption}. Indeed, the terminal utility depends on the net wealth after accounting for the liability, suggesting that the continuation value should incorporate a process representing the discounted future liability exposure. Motivated by this observation, and following the BSDE approach to (recursive) utility maximisation, we introduce an auxiliary process $Y$ and seek an optimal value process of the form:
\begin{align}\label{Ansatz_Optimal process}
	\widehat{V}_{t}:=e^{-\delta\theta t}\frac{(X_{t}^{c,\pi}+e^{\int_{0}^{t}r_{s}\mathrm{d}s}Y_{t}) ^{1-\gamma}}{1-\gamma}~\text{ for }~t\in[0,T],
\end{align}
where $Y$ will be chosen so that the terminal condition $X_{T}^{c,\pi}+e^{\int_{0}^{T}r_{s}\mathrm{d}s}Y_{T}=X_{T}^{c,\pi}-K_{T}$
is satisfied and $X_{0}^{c,\pi}+Y_{0}$ ($=x+Y_{0}$) is a constant independent of $(c,\pi)$. We therefore define $Y$ as the first part of the solution of the BSDE
\begin{align}\label{Auxiliary BSDE}
	\mathrm{d}Y_{t}=-\mathcal{H}(t,X_{t}^{c,\pi},Y_{t},Z_{t})\mathrm{d}t+Z_{t}\mathrm{d}I_{t},~~Y_{T}=-K_{T}e^{-\int_{0}^{T}r_{s}\mathrm{d}s}
\end{align}
for some generator $\mathcal{H}$ to be defined. In line with the permissible sets considered in, for example, \cite{aurand2023epstein,feng2026consumption,xing2017consumption}, the set of permissible consumption-portfolio strategies is defined as follows.
\begin{defi}\label{Permissible strategies}
	An admissible consumption-portfolio strategy $(c,\pi)$ of is \textit{permissible} if
	\begin{itemize}
		\item[$(i)$] $c\in\mathcal{C}_{a}$ with $c_{T}=X_{T}^{c,\pi}+e^{\int_{0}^{T}r_{s}\mathrm{d}s}Y_{T}$;
		\item[$(ii)$] $X_{t}^{c,\pi}+e^{\int_{0}^{t}r_{s}\mathrm{d}s}Y_{t}>0$ for all $t\in[0,T]$;
		\item[$(iii)$] $(X_{\cdot}^{c,\pi}+e^{\int_{0}^{\cdot}r_{s}\mathrm{d}s}Y_{\cdot})^{1-\gamma}$ is of class (D) on $[0,T]$.
	\end{itemize}
\end{defi}
We denote by $\mathcal{P}$ the set of permissible consumption--portfolio strategies. The following remark provides further justification for the permissibility conditions defining $\mathcal{P}$.
\begin{remark}
	The permissibility conditions in Definition \ref{Permissible strategies} each play a distinct role. Condition $(i)$ is standard in the literature on Epstein--Zin recursive utility. It guarantees the well-posedness of the utility process and requires that the investor's total wealth, net of the terminal liability, is fully consumed at the terminal time. This terminal consumption condition is essential for establishing an explicit connection between the utility process and the investment strategy. Condition $(ii)$ ensures the well-posedness of the ansatz \eqref{Ansatz_Optimal process}, while Condition $(iii)$ is indispensable for verifying the optimality of the candidate strategy. For a more detailed discussion of the role of the class (D) property, we refer the reader to \cite[Rmk.~3.6]{aurand2023epstein}.
\end{remark}

From the martingale optimality principle (see, e.g., \cite{hu2005utility,xing2017consumption}), we must choose the function $\mathcal{H}$ in \eqref{Auxiliary BSDE} such that the process
\begin{align}\label{Optimal martingale}
	M_{t}^{c,\pi}:=e^{-\delta\theta t}\frac{(X_{t}^{c,\pi}+e^{\int_{0}^{t}r_{s}\mathrm{d}s}Y_{t})^{1-\gamma}}{1-\gamma}+\int_{0}^{t}f\big(c_{s},{e^{-\delta\theta s}}\frac{(X_{s}^{c,\pi}+e^{\int_{0}^{s}r_{u}\mathrm{d}u}Y_{s})^{1-\gamma}}{1-\gamma}\big)\mathrm{d}s
\end{align}
for $t\in[0,T]$, is a local supermartingale for all $(c,\pi)\in\mathcal{P}$ and there exists $(c^{*},\pi^{*})\in\mathcal{P}$ such that $M^{c^{*},\pi^{*}}$ is a local martingale. It{\^o}'s formula applied to $M^{c,\pi}$ gives
\begin{align}\label{OptEZ_Differential_Insurer}
	\mathrm{d}M_{t}^{c,\pi}&=e^{-\delta\theta t} (X_{t}^{c,\pi}+e^{\int_{0}^{t}r_{s}\mathrm{d}s}Y_{t})^{-\gamma}\Big(-c_{t}+\delta\frac{c_{t}^{1-\frac{1}{\psi}}}{1-\frac{1}{\psi}}(X_{t}^{c,\pi}+e^{\int_{0}^{t}r_{s}\mathrm{d}s}Y_{t})^{\frac{1}{\psi}}-e^{\int_{0}^{t}r_{s}\mathrm{d}s}Z_{t}^{\top}\eta_{t}\notag\\
	&+\frac{1}{2\gamma}(X_{t}^{c,\pi}+e^{\int_{0}^{t}r_{s}\mathrm{d}s}Y_{t})\|\eta_{t}\|^{2}+r_{t}(X_{t}^{c,\pi}+e^{\int_{0}^{t}r_{s}\mathrm{d}s}Y_{t})-\frac{\delta\theta}{1-\gamma}(X_{t}^{c,\pi}+e^{\int_{0}^{t}r_{s}\mathrm{d}s}Y)\notag\\
	&-e^{\int_{0}^{t}r_{s}\mathrm{d}s}\mathcal{H}(t,X_{t}^{c,\pi},Y_{t},Z_{t})\Big)\mathrm{d}t\notag\\
	&-\frac{\gamma}{2}e^{-\delta\theta t} (X_{t}^{c,\pi}+e^{\int_{0}^{t}r_{s}\mathrm{d}s}Y_{t})^{-\gamma-1}\Big\|\pi_{t}+\Big(e^{\int_{0}^{t}r_{s}\mathrm{d}s}Z_{t}-\frac{1}{\gamma}(X_{t}^{c,\pi}+e^{\int_{0}^{t}r_{s}\mathrm{d}s}Y_{t})\eta_{t}\Big)\Big\|^{2}\mathrm{d}t\notag\\
	&+e^{-\delta\theta t} (X_{t}^{c,\pi}+e^{\int_{0}^{t}r_{s}\mathrm{d}s}Y_{t})^{-\gamma}\big(\pi_{t}^{\top}+e^{\int_{0}^{t}r_{s}\mathrm{d}s}Z_{t}^{\top}\big)\mathrm{d}I_{t}.
\end{align}
Expecting the drift to be non-positive for any $(c,\pi)\in\mathcal{P}$ and zero at an optimal strategy $(c^{*},\pi^{*})\in\mathcal{P}$, the generator $\mathcal{H}$ can be obtained by formally taking the supremum over $\pi$ and $c$ in the drift of \eqref{OptEZ_Differential_Insurer} and setting it to be zero. Hence, a candidate optimal strategy is given by
\begin{align}\label{Candidate optimal strategy}
	\pi_{t}^{*}=-e^{\int_{0}^{t}r_{s}\mathrm{d}s}Z_{t}+\frac{1}{\gamma}(X_{t}^{c^{*},\pi^{*}}+e^{\int_{0}^{t}r_{s}\mathrm{d}s}Y_{t})\eta_{t},~~c_{t}^{*}=\delta^{\psi}(X_{t}^{c^{*},\pi^{*}}+e^{\int_{0}^{t}r_{s}\mathrm{d}s}Y_{t}),~0\le t<T,
\end{align}
where $(X^{c^{*},\pi^{*}},Y,Z)$ is solution to the coupled linear FBSDE
\begin{align}\label{Optimal FBSDE}
	\begin{cases}
		\mathrm{d}Y_{t}&=-\Big(e^{-\int_{0}^{t}r_{s}\mathrm{d}s}\Big(\frac{\delta^{\psi}}{\psi-1}+r_{t}+\frac{1}{2\gamma}\|\eta_{t}\|^{2}{-\frac{\delta\theta}{1-\gamma}}\Big)(X_{t}^{c^{*},\pi^{*}}+e^{\int_{0}^{t}r_{s}\mathrm{d}s}Y_{t})\\
		&\phantom{XXx}-Z_{t}^{\top}\eta_{t}\Big)\mathrm{d}t+Z_{t}\mathrm{d}I_{t}\\
		\mathrm{d}X_{t}^{c^{*},\pi^{*}}&=\Big(r_{t}X_{t}^{c^{*},\pi^{*}}+\big(-\delta^{\psi}+\frac{1}{\gamma}\|\eta_{t}\|^{2}\big)(X_{t}^{c^{*},\pi^{*}}+e^{\int_{0}^{t}r_{s}\mathrm{d}s}Y_{t})-e^{\int_{0}^{t}r_{s}\mathrm{d}s}Z_{t}^{\top}\eta_{t}\Big)\mathrm{d}t\\
		&\phantom{x}+\Big(\frac{1}{\gamma}(X_{t}^{c^{*},\pi^{*}}+e^{\int_{0}^{t}r_{s}\mathrm{d}s}Y_{t})\eta_{t}^{\top}-e^{\int_{0}^{t}r_{s}\mathrm{d}s}Z_{t}^{\top}\Big)\mathrm{d}I_{t}\\
		X_{0}^{c^{*},\pi^{*}}&=x\ge0,~~Y_{T}=-K_{T}e^{-\int_{0}^{T}r_{s}\mathrm{d}s}.
	\end{cases}
\end{align}

Therefore, the candidate solution to problem \eqref{Main problem} is given by \eqref{Candidate optimal strategy}, provided that the coupled FBSDE \eqref{Optimal FBSDE} with unbounded random coefficients is well-defined in an appropriate function space.
\begin{remark}\label{Some remarks on the solvability of FBSDEs}
	We use a decoupling reduction method via which an explicit solution to the coupled linear FBSDE with unbounded coefficients associated to the optimal consumption and investment allocations can be obtained. This is due to the special structure of the FBSDE. In contrast to existing works such as \cite{ma1994solving,yong1999linear,yong2006linear} in which solvability to a coupled linear FBSDE with bounded coefficients is reduced to the solvability of BSDEs, in the present paper we instead reduce the solvability of our coupled linear FBSDE with unbounded coefficients to that of Forward SDEs. Moreover, the unboundedness of the coefficients precludes the use of \cite{ankirchner2022transformation,hu1995solution,ma2015well,pardoux1990adapted,pardoux1999forward,peng1999fully,xie2020exploration}.
\end{remark}

To show the well-definedness and construct a solution of the FBSDE \eqref{Optimal FBSDE}, we consider the following conditions.
\begin{assum}\leavevmode\label{Change of measure}
	\begin{itemize}
		\item[$(i)$] For $q>1$, define $q_{max}:=\max\Big\{\frac{45+5\gamma}{\gamma^{2}},\frac{2q^{2}+q(\gamma-1)}{\gamma^{2}},\frac{9(2q-\gamma)^{2}}{2\gamma^{2}},3\frac{\gamma^{2}-5q\gamma+6q^{2}-q}{\gamma^{2}},32\Big\}$. Suppose the market price of risk $\eta$ satisfies
		\begin{align}
			\mathbb{E}\Big[\exp\Big(q_{max}\int_{0}^{T}\|\eta_{s}\|^{2}\mathrm{d}s\Big)~\big|\mathcal{F}_{0}^{r,S}\Big]&<\infty.
		\end{align}
		\item[$(ii)$]  $K_{T}e^{-\int_{0}^{T}r_{s}\mathrm{d}s}\in\mathcal{L}_{\mathbb{Q}_{0}}^{2q},~{q>1}$, where $\mathbb{Q}_{0}$ is the probability measure equivalent to $\mathbb{P}$ defined by $\frac{\mathrm{d}\mathbb{Q}_{0}}{\mathrm{d}\mathbb{P}}\big|_{\mathcal{F}_{T}^{r,S}}:=\mathcal{E}\big(\int-\eta^{\top}\mathrm{d}I\big)_{T}:=\exp\left(-\frac{1}{2}\int_{0}^{T}\|\eta_{s}\|^{2}\mathrm{d}s-\int_{0}^{T}\eta_{s}^{\top}\mathrm{d}I_{s} \right)$.
	\end{itemize}
\end{assum}
\begin{remark}\label{Some remarks on the exponential integrability conditions}
	Imposing exponential integrability conditions to establish existence results of (F)BSDEs with unbounded coefficients is also considered by many authors; see, for example, assumption $2.1$ in \cite{feng2026consumption,sun2020mean,xing2017consumption}.
\end{remark}

We define the processes $(H_{t})_{t\in[0,T]}$, $(\alpha_{t})_{t\in[0,T]}$ and $(\varphi_{t})_{t\in[0,T]}$ by
\begin{align}\label{Auxiliary processes}
	\begin{cases}
		&H_{t}:=\mathcal{E}\big(\int-\eta^{\top}\mathrm{d}I\big)_{t},\quad\alpha_{t}:=e^{-\int_{0}^{t}r_{s}\mathrm{d}s}\Big(\frac{\delta^{\psi}}{\psi-1}+r_{t}+\frac{1}{2\gamma}\|\eta_{t}\|^{2}-\frac{\delta\theta}{1-\gamma}\Big)\\
		&\text{and }\varphi_{t}:=\exp\Big(\int_{0}^{t}\Big(-\frac{\delta^{\psi}\psi} {\psi-1}+\frac{\gamma-1}{2\gamma^{2}}\|\eta_{s}\|^{2}+\frac{\delta\theta}{1-\gamma}\Big)\mathrm{d}s+\frac{1}{\gamma}\int_{0}^{t}\eta_{s}^{\top}\mathrm{d}I_{s}\Big).
	\end{cases}
\end{align}
\begin{remark}\label{Integrability of alphaphi}
	Assumption~\ref{Change of measure}$.(i)$ yields $\alpha\varphi\in\mathcal{H}_{\mathbb{Q}_{0}}^{2q},~{q>1}$, (see Appendix~\ref{Appendix A}). This is used in the construction of a solution to the FBSDE \eqref{Optimal FBSDE}; see Proposition~\ref{Existence result for the FBSDE}.
\end{remark}
\begin{prop}\label{Existence result for the FBSDE}
	Suppose Assumption \ref{Change of measure} holds. Let $\widetilde{x}$ denotes the constant defined by $\widetilde{x}:=\frac{x-\mathbb{E}\big[H_{T}K_{T}e^{-\int_{0}^{T}r_{s}\mathrm{d}s}~|\mathcal{F}_{0}^{r,S}\big]}{1-\mathbb{E}\big[\int_{0}^{T}H_{s}\alpha_{s}\varphi_{s}\mathrm{d}s~|\mathcal{F}_{0}^{r,S}\big]}$. Then the FBSDE \eqref{Optimal FBSDE} admits a solution\\ $(X^{c^{*},\pi^{*}},Y,Z)\in\mathcal{H}_{\mathbb{P}}^{q'}\times{\mathcal{S}_{\mathbb{P}}^{q}}\times\mathbb{H}_{\mathbb{P}}^{q},~q>q'>1$, satisfying
	\begin{align}\label{Relation between X and Y}
		X_{t}^{c^{*},\pi^{*}}+e^{\int_{0}^{t}r_{s}\mathrm{d}s}Y_{t}=\widetilde{x}\varphi_{t},~0\le t\le T.
	\end{align}
	Moreover, $(Y,Z)\in{\mathcal{S}_{\mathbb{P}}^{q}}\times\mathbb{H}_{\mathbb{P}}^{q},~{q>1}$, is the unique solution to the BSDE
	\begin{align}\label{Auxiliary BSDE_Investor}
		\mathrm{d}Y_{t}&=-\Big(\widetilde{x}\alpha_{t}\varphi_{t}-Z_{t}^{\top}\eta_{t}\Big)\mathrm{d}t+Z_{t}^{\top}\mathrm{d}I_{t},~~Y_{T}=-K_{T}e^{-\int_{0}^{T}r_{s}\mathrm{d}s},
	\end{align}
	In addition, the expectation representation of the first component $Y$ is given by
	\begin{align}\label{YSolution to the BSDE}
		Y_{t}&=H_{t}^{-1}\mathbb{E}\Big[-H_{T}K_{T}e^{-\int_{0}^{T}r_{s}\mathrm{d}s}+\widetilde{x}\int_{t}^{T}H_{s}\alpha_{s}\varphi_{s}\mathrm{d}s~\big|\mathcal{F}_{t}^{r,S}\Big],\,\,0\le t\le T.
	\end{align}
\end{prop}
\begin{proof}
	Due to the length of the proof, we defer it to Appendix \ref{Appendix B}.
\end{proof}

We are now ready to give the main result of this paper
\begin{thm}\label{Mainresult1}
	Assume $x>\mathbb{E}\big[H_{T}K_{T}e^{-\int_{0}^{T}r_{s}\mathrm{d}s}~\big|\mathcal{F}_{0}^{r,S}\big]$ and Assumption \ref{Change of measure} holds. Let $\widetilde{x}$ be the constant defined in Proposition~\ref{Existence result for the FBSDE}. Then the optimal consumption and portfolio strategy for the optimisation problem \eqref{Main problem} is given by
	\begin{align}\label{Optimal strategy}
		c_{t}^{*}=\delta^{\psi}\widetilde{x}\varphi_{t}~\text{ and }~\pi_{t}^{*}=-e^{\int_{0}^{t}r_{s}\mathrm{d}s}Z_{t}+\frac{1}{\gamma}\widetilde{x}\varphi_{t}\eta_{t}.
	\end{align}
	In particular, the corresponding holdings in the stocks and the bond are given by\\ $(\pi_{t}^{S},\pi_{t}^{B})=(\pi_{t}^{*})^{\top}\Sigma_{t}^{-1}$ for $t\in[0,T]$, where $\Sigma$ is defined below \eqref{Wealth}.
	
	Moreover, the value function of problem \eqref{Main problem} is given by 
	\begin{align}\label{Optimal value function}
		\widehat{V}_{0}&=\frac{1}{1-\gamma}\left(\frac{x-\mathbb{E}\big[H_{T}K_{T}e^{-\int_{0}^{T}r_{s}\mathrm{d}s}~\big|\mathcal{F}_{0}^{r,S}\big]}{1-\mathbb{E}\big[\int_{0}^{T}H_{s}\alpha_{s}\varphi_{s}\mathrm{d}s~\big|\mathcal{F}_{0}^{r,S}\big]}\right)^{1-\gamma}.
	\end{align}
\end{thm}
\begin{proof}
	First, we prove that $(c^{*},\pi^{*})\in\mathcal{P}$. (Recall $\mathcal{P}$ from Definition \ref{Permissible strategies}). Clearly, $X_{t}^{c^{*},\pi^{*}}+e^{\int_{0}^{t}r_{s}\mathrm{d}s}Y_{t}=\widetilde{x}\varphi_{t}>0,~t\in[0,T]$; due to $x>\mathbb{E}\big[H_{T}K_{T}e^{-\int_{0}^{T}r_{s}\mathrm{d}s}~\big|\mathcal{F}_{0}^{r,S}\big]$ and \eqref{E1_Proposition 3.2}. Besides,
	\begin{align}\label{Class (D) property}
		\big(X_{t}^{c^{*},\pi^{*}}+e^{\int_{0}^{t}r_{s}\mathrm{d}s}Y_{t}\big)^{1-\gamma}&=\widetilde{x}^{1-\gamma}\exp\Big(\int_{0}^{t}\big(-\delta^{\psi}\theta e^{-\delta\theta\psi s}+\delta\theta\big)\mathrm{d}s\Big)\mathcal{E} \big(\int\frac{1-\gamma}{\gamma}\eta^{\top}\mathrm{d}I\big)_{t}.
	\end{align}
	Using Assumption~\ref{Change of measure}$.(i)$ with $\big(\frac{1-\gamma}{\gamma}\big)^{2}<1<q_{max}$, we deduce that $\mathcal{E}\big(\int\frac{1-\gamma}{\gamma}\eta^{\top}\mathrm{d}I\big)$ is a $\mathbb{P}$-martingale (hence of class (D)). Thus the right-side of \eqref{Class (D) property} is of class (D) as a product of a bounded deterministic function and a process of class (D). Therefore, $(X^{c^{*},\pi^{*}}+e^{\int_{0}^{}r_{s}\mathrm{d}s}Y)^{1-\gamma}$ is of class (D) on $[0,T]$. Finally, using \cite[Prop.~2.2]{xing2017consumption} and the latter class (D) property, to show that $c\in\mathcal{C}_{a}$ it suffices to prove that $\mathbb{E}\big[\int_{0}^{T}(X_{t}^{c^{*},\pi^{*}}+e^{\int_{0}^{t}r_{s}\mathrm{d}s}Y_{t})^{1-\frac{1}{\psi}}\mathrm{d}t~\big|\mathcal{F}_{0}^{r,S}\big]<\infty$. If $\gamma\psi=1,\gamma>1$, then the latter inequality follows from \eqref{Class (D) property}. If $\gamma>1,\psi>1$, then using successively Cauchy-Schwarz inequality, the inequality $\exp\big(\int_{0}^{t}\big(-\delta^{\psi}+\frac{\delta}{1-\gamma}\big)\mathrm{d}s\big)\le\exp\big(\big|\frac{\delta}{1-\gamma}\big|T\big)$ for $t\in[0,T]$, Novikov condition and Assumption~\ref{Change of measure}$.(i)$ with $\big(1-\frac{1}{\psi}\big)\big(\frac{\gamma+1}{\gamma}-\frac{2}{\gamma\psi^{2}}\big)<2<q_{max}$ and $\big(1-\frac{1}{\psi}\big)^{2}\frac{2}{\gamma^{2}}<2<q_{max}$, we obtain
	\begin{align*}
		&\mathbb{E}\Big[\int_{0}^{T}\big(X_{t}^{c^{*},\pi^{*}}+e^{\int_{0}^{t}r_{s}\mathrm{d}s}Y_{t}\big)^{1-\frac{1}{\psi}}\mathrm{d}t~\big|\mathcal{F}_{0}^{r,S}\Big]\notag\\
		&\le\Big(\mathbb{E}\Big[\int_{0}^{T}\exp\Big(\big(1-\frac{1}{\psi}\big)\big(\frac{\gamma+1}{\gamma}-\frac{2}{\gamma\psi^{2}}\big)\int_{0}^{t}\|\eta_{s}\|^{2}\mathrm{d}s\Big)\mathrm{d}t~\big|\mathcal{F}_{0}^{r,S}\Big]\Big)^{\frac{1}{2}}\notag\\
		&\phantom{xx}\times\Big(\mathbb{E}\Big[\int_{0}^{T}\mathcal{E}\big(\int\big(1-\frac{1}{\psi}\big)\frac{2}{\gamma}\eta^{\top}\mathrm{d}I\big)_{t}\mathrm{d}t~\big|\mathcal{F}_{0}^{r,S}\Big]\Big)^{\frac{1}{2}}\exp\big(\big|\frac{\delta}{1-\gamma}\big|T\big)\widetilde{x}^{1-\frac{1}{\psi}}<\infty.
	\end{align*}
	
	Second, we show that $(c^{*},\pi^{*})$ is optimal. The proof follows similar arguments as in the proof of proposition $3.2$ in \cite{feng2026consumption}. Recall that $M^{c,\pi}$ is a local supermartingale for all $(c,\pi)\in\mathcal{P}$ and a local martingale for $(c^{*},\pi^{*})$. Then, for a non-decreasing sequence $(\tau_{n})_{n\geq1}$ of stopping times such that $\tau_{n}\stackrel{a.s.}\longrightarrow T$ and for all $(c,\pi)\in\mathcal{P}$, we have
	\begin{align*}
		&\mathbb{E}\left[M_{\tau_{n}}^{c,\pi}\right]\notag\\
		&=\mathbb{E}\left[e^{-\delta\theta\tau_{n}}\frac{(X_{\tau_{n}}^{c,\pi}+e^{\int_{0}^{\tau_{n}}r_{s}\mathrm{d}s}Y_{\tau_{n}})^{1-\gamma}}{1-\gamma}+\int_{0}^{\tau_{n}}f\Big(c_{s},e^{-\delta\theta s} \frac{(X_{s}^{c,\pi}+e^{\int_{0}^{s}r_{u}\mathrm{d}u}Y_{s}) ^{1-\gamma}}{1-\gamma}\Big)\mathrm{d}s\right]\\
		&\leq M_{0}^{c,\pi}\notag\\
		&=M_{0}^{c^{*},\pi^{*}}\notag\\
		&=\mathbb{E}\left[e^{-\delta\theta\tau_{n}}\frac{(X_{\tau_{n}}^{c^{*},\pi^{*}}+e^{\int_{0}^{\tau_{n}}r_{s}\mathrm{d}s}Y_{\tau_{n}})^{1-\gamma}}{1-\gamma}+\int_{0}^{\tau_{n}}f\Big(c_{s}^{*},e^{-\delta\theta s} \frac{(X_{s}^{c^{*},\pi^{*}}+e^{\int_{0}^{s}r_{u}\mathrm{d}u}Y_{s})^{1-\gamma}}{1-\gamma}\Big)\mathrm{d}s\right].
	\end{align*}
	Hence,
	\begin{align*}
		&\mathbb{E}\left[e^{-\delta\theta\tau_{n}}\frac{(X_{\tau_{n}}^{c,\pi}+e^{\int_{0}^{\tau_{n}}r_{s}\mathrm{d}s}Y_{\tau_{n}})^{1-\gamma}}{1-\gamma}+\int_{0}^{\tau_{n}}f\Big(c_{s},e^{-\delta\theta s} \frac{(X_{s}^{c,\pi}+e^{\int_{0}^{s}r_{u}\mathrm{d}u}Y_{s}) ^{1-\gamma}}{1-\gamma}\Big)\mathrm{d}s\right]\\
		&\leq\mathbb{E}\left[e^{-\delta\theta\tau_{n}}\frac{(X_{\tau_{n}}^{c^{*},\pi^{*}}+e^{\int_{0}^{\tau_{n}}r_{s}\mathrm{d}s}Y_{\tau_{n}})^{1-\gamma}}{1-\gamma}+\int_{0}^{\tau_{n}}f\Big(c_{s}^{*},e^{-\delta\theta s} \frac{(X_{s}^{c^{*},\pi^{*}}+e^{\int_{0}^{s}r_{u}\mathrm{d}u}Y_{s})^{1-\gamma}}{1-\gamma}\Big)\mathrm{d}s\right].
	\end{align*}
	Since $(c^*,\pi^*)\in\mathcal{P}$, sending $n\to\infty$, the monotone convergence theorem and the class (D) property of $(X_{\cdot}^{c,\pi}+e^{\int_{0}^{\cdot}r_{s}\mathrm{d}s}Y_{\cdot})^{1-\gamma}$ and $(X_{\cdot}^{c^{*},\pi^{*}}+e^{\int_{0}^{\cdot}r_{s}\mathrm{d}s}Y_{\cdot})^{1-\gamma}$ yields
	\begin{align}\label{maxuti}
		&\sup_{(c,\pi)\in\mathcal{P}}\mathbb{E}\left[e^{-\delta\theta T}\frac{(X_{T}^{c,\pi}+e^{\int_{0}^{T}r_{s}\mathrm{d}s}Y_{T})^{1-\gamma}}{1-\gamma}+\int_{0}^{T}f\Big(c_{s},e^{-\delta\theta s} \frac{(X_{s}^{c,\pi}+e^{\int_{0}^{s}r_{u}\mathrm{d}u}Y_{s}) ^{1-\gamma}}{1-\gamma}\Big)\mathrm{d}s\right]\notag\\
		&\leq\mathbb{E}\left[e^{-\delta\theta T}\frac{(X_{T}^{c^{*},\pi^{*}}+e^{\int_{0}^{T}r_{s}\mathrm{d}s}Y_{T})^{1-\gamma}}{1-\gamma}+\int_{0}^{T}f\Big(c_{s}^{*},e^{-\delta\theta s} \frac{(X_{s}^{c^{*},\pi^{*}}+e^{\int_{0}^{s}r_{u}\mathrm{d}u}Y_{s})^{1-\gamma}}{1-\gamma}\Big)\mathrm{d}s\right].
	\end{align}
	The arguments leading to \eqref{maxuti} is valid in both the empirically relevant regime $\gamma>1,\,\psi>1$ and the CRRA case, characterised by $\gamma>1$ and $\gamma\psi=1$. Thus, the pair $(c^{*},\pi^{*})$ is optimal for the optimisation problem \eqref{Main problem}.
\end{proof}

Our next objective is to establish that the solution to the BSDE \eqref{Auxiliary BSDE_Investor} is Malliavin differentiable; see \cite{alos2024malliavin,nunno2008malliavin,nualart2006malliavin} for comprehensive introductions on Malliavin calculus. This regularity property allows the optimal investment strategy $\pi^{*}$, given in \eqref{Optimal strategy}, to be expressed in term of conditional expectations of random variables that can be efficiently approximated via Monte Carlo simulation. In the sequel, $D_{t}(\cdot)$ denotes the Malliavin operator at all $t\in[0,T]$. We work under the following conditions.
\begin{assum}\label{Conditions to be Malliavin differentiable}\leavevmode
	\begin{itemize}
		\item[$(i)$] $K_{T}e^{-\int_{0}^{T}r_{s}\mathrm{d}s}\in\mathbb{D}^{1,2}$, $H_{T}K_{T}e^{-\int_{0}^{T}r_{s}\mathrm{d}s}\in\mathbb{D}^{1,2}$.
		\item[$(ii)$] $\mathbb{E}^{\mathbb{Q}_{0}}\Big[\big|K_{T}\big|e^{-\int_{0}^{T}r_{s}\mathrm{d}s}~\big|\mathcal{F}_{0}^{r,S}\Big]<\infty$.
		\item[$(iii)$] $\mathbb{E}^{\mathbb{Q}_{0}}\Big[\int_{0}^{T}\big\|D_{t}\big(K_{T}e^{-\int_{0}^{T}r_{s}\mathrm{d}s}\big)\big\|^{2}\mathrm{d}t~\big|\mathcal{F}_{0}^{r,S}\Big]<\infty$.
	\end{itemize}
\end{assum}
\begin{prop}\label{Auxilliary results_Malliavin differentiability}
	Let Assumptions \ref{Change of measure} and \ref{Conditions to be Malliavin differentiable} hold. Then, the following results hold:
	\begin{itemize}
		\item[$(i)$] $K_{T}e^{-\int_{0}^{T}r_{s}\mathrm{d}s}+\int_{0}^{T}\alpha_{s}\varphi_{s}\mathrm{d}s\in\mathbb{D}^{1,2}$, $H_{T}\Big(K_{T}e^{-\int_{0}^{T}r_{s}\mathrm{d}s}+\int_{0}^{T}\alpha_{s} \varphi_{s}\mathrm{d}s\Big)\in\mathbb{D}^{1,2}$.
		\item[$(ii)$] $\eta_{t}\in\mathbb{D}^{1,2}$ for almost all $t\in[0,T]$.
		\item[$(iii)$] $\mathbb{E}^{\mathbb{Q}_{0}}\Big[\big|K_{T}e^{-\int_{0}^{T}r_{s}\mathrm{d}s}+\int_{0}^{T}\alpha_{s}\varphi_{s}\mathrm{d}s\big|~\big|\mathcal{F}_{0}^{r,S}\Big]<\infty$.
		\item[$(iv)$] $\mathbb{E}^{\mathbb{Q}_{0}}\Big[\int_{0}^{T}\Big(\big\|D_{t}\big(K_{T}e^{-\int_{0}^{T}r_{s}\mathrm{d}s}\big)\big\|^{2}+\big\|D_{t}\big(\int_{0}^{T}\alpha_{s}\varphi_{s}\mathrm{d}s\big)\big\|^{2}\Big)\mathrm{d}t~\big|\mathcal{F}_{0}^{r,S}\Big]<\infty$.
		\item[$(v)$] $\big(D_{t}(\alpha_{t}\varphi_{t})-Z_{t}^{\top}D_{t}(\eta_{t})\big)_{t\in[0,T]}\in\mathbb{H}_{\mathbb{Q}_{0}}^{2}$.
	\end{itemize}
\end{prop}
\begin{proof}
	Due to the length of the proof, we postpone it in Appendix \ref{Appendix C}.
\end{proof}

Results $(i)$-$(iv)$ in Proposition \ref{Auxilliary results_Malliavin differentiability} ensure the applicability of the Clark-Ocone formula to the $\mathcal{F}_{T}^{r,S}$-random variable $K_{T}e^{-\int_{0}^{T}r_{s}\mathrm{d}s}+\widetilde{x}\int_{0}^{T}\alpha_{s}\varphi_{s}\mathrm{d}s$ under the new measure $\mathbb{Q}_{0}$. Compare with \cite[Thm.~4.5, Rmk.~4.6]{nunno2008malliavin}.
\begin{prop}\label{ZSolution to the auxiliary BSDE_2}
	Let Assumptions \ref{Change of measure} and \ref{Conditions to be Malliavin differentiable} hold. Then, the unique solution\\ $(Y,Z)\in{\mathcal{S}_{\mathbb{P}}^{q}}\times\mathbb{H}_{\mathbb{P}}^{q},~{q>1}$, to the BSDE \eqref{Auxiliary BSDE_Investor} is Malliavin differentiable and we have
	\begin{align}\label{ZSolution to the FBSDE1}
		Z_{t}&=D_{t}(Y_{t}),~\text{where $D_{t}(\cdot)$ denote the Malliavin operator for all $t\in[0,T]$.}
	\end{align}
\end{prop}
\begin{proof}
	We define the processes $\tilde{Y}_{t}:=Y_{t}+\widetilde{x}\int_{0}^{t}\alpha_{s}\varphi_{s}\mathrm{d}s$ and $\tilde{Z}_{t}:=Z_{t}$ for $t\in[0,T]$. Hence, $(\tilde{Y},\tilde{Z})$ is the unique solution to the BSDE
	\begin{align}\label{E3_Proposition 3.9}
		\mathrm{d}\tilde{Y}_{t}&=\tilde{Z}_{t}^{\top}\eta_{t}\mathrm{d}t+\tilde{Z}_{t}^{\top}\mathrm{d}I_{t}=\tilde{Z}_{t}^{\top}\mathrm{d}I_{t}^{\mathbb{Q}_{0}},~~\tilde{Y}_{T}=-K_{T}e^{-\int_{0}^{T}r_{s}\mathrm{d}s}+\widetilde{x}\int_{0}^{T}\alpha_{s}\varphi_{s}\mathrm{d}s,
	\end{align}
	where $I_{t}^{\mathbb{Q}_{0}}:=I_{t}+\int_{0}^{t}\eta_{s}\mathrm{d}s,\,t\in[0,T]$, is a Brownian motion under $\mathbb{Q}_{0}$. Then
	\begin{align}\label{E4_Proposition 3.9}
		-K_{T}e^{-\int_{0}^{T}r_{s}\mathrm{d}s}+\widetilde{x}\int_{0}^{T}\alpha_{s}\varphi_{s}\mathrm{d}s&=Y_{0}+\int_{0}^{T}\tilde{Z}_{s}^{\top}\mathrm{d}I_{s}^{\mathbb{Q}_{0}}.
	\end{align}
	Using Proposition \ref{Auxilliary results_Malliavin differentiability} and applying the Clark-Ocone formula under change of measure as in \cite[Thm.~4.5]{nunno2008malliavin} to $\tilde{Y}_{T}=-K_{T}e^{-\int_{0}^{T}r_{s}\mathrm{d}s}+\widetilde{x}\int_{0}^{T}\alpha_{s}\varphi_{s}\mathrm{d}s\in\mathbb{D}^{1,2}$, we obtain
	\begin{align}\label{E5_Proposition 3.9}
		&-K_{T}e^{-\int_{0}^{T}r_{s}\mathrm{d}s}+\widetilde{x}\int_{0}^{T}\alpha_{s}\varphi_{s}\mathrm{d}s\notag\\
		&=\mathbb{E}^{\mathbb{Q}_{0}}\Big[-K_{T}e^{-\int_{0}^{T}r_{s}\mathrm{d}s}+\widetilde{x}\int_{0}^{T}\alpha_{s}\varphi_{s}\mathrm{d}s~\big|\mathcal{F}_{0}^{r,S}\Big]\notag\\
		&\phantom{xx}+\int_{0}^{T}\mathbb{E}^{\mathbb{Q}_{0}}\Big[D_{t}\big(-K_{T}e^{-\int_{0}^{T}r_{s}\mathrm{d}s}+\widetilde{x}\int_{0}^{T}\alpha_{s}\varphi_{s}\mathrm{d}s\big)\notag\\
		&\phantom{xx}-\big(-K_{T}e^{-\int_{0}^{T}r_{s}\mathrm{d}s}+\widetilde{x}\int_{0}^{T}\alpha_{s}\varphi_{s}\mathrm{d}s\big)\int_{t}^{T}D_{t}(\eta_{s})\mathrm{d}I_{s}^{\mathbb{Q}_{0}}~|\mathcal{F}_{t}^{r,S}\Big]^{\top}\mathrm{d}I_{t}^{\mathbb{Q}_{0}}.
	\end{align}
	By uniqueness of the solution to the BSDE \eqref{E3_Proposition 3.9}, we deduce from \eqref{E4_Proposition 3.9}-\eqref{E5_Proposition 3.9} that
	\begin{align}\label{E6_Proposition 3.9}
		Y_{0}&=\mathbb{E}\Big[-H_{T}K_{T}e^{-\int_{0}^{T}r_{s}\mathrm{d}s}+\widetilde{x}\int_{0}^{T}H_{s}\alpha_{s}\varphi_{s}\mathrm{d}s~\big|\mathcal{F}_{0}^{r,S}\Big]
	\end{align}
	as we already obtained in Proposition~\ref{Existence result for the FBSDE}, and
	\begin{align}\label{E8_Proposition 3.9}
		Z_{t}=\tilde{Z}_{t}&=\mathbb{E}^{\mathbb{Q}_{0}}\Big[D_{t}\big(-K_{T}e^{-\int_{0}^{T}r_{s}\mathrm{d}s}+\widetilde{x}\int_{0}^{T}\alpha_{s}\varphi_{s}\mathrm{d}s\big)\notag\\
		&\phantom{XXXXx}-\big(-K_{T}e^{-\int_{0}^{T}r_{s}\mathrm{d}s}+\widetilde{x}\int_{0}^{T}\alpha_{s}\varphi_{s}\mathrm{d}s\big)\int_{t}^{T}D_{t}(\eta_{s})\mathrm{d}I_{s}^{\mathbb{Q}_{0}}~|\mathcal{F}_{t}^{r,S}\Big].
	\end{align}
	Besides, we consider the BSDE 
	\begin{align}\label{E9_Proposition 3.9}
		\begin{cases}
			\mathrm{d}D_{t}(Y_{t})&=-\big(\widetilde{x}D_{t}(\alpha_{t}\varphi_{t})-D_{t}(Z_{t}^{\top})\eta_{t}-Z_{t}^{\top}D_{t}(\eta_{t})\big)\mathrm{d}t+D_{t}(Z_{t}^{\top})\mathrm{d}I_{t}\\
			D_{t}(Y_{T})&=D_{t}(-K_{T}e^{-\int_{0}^{T}r_{s}\mathrm{d}s}).
		\end{cases}
	\end{align}
	Using Proposition \ref{Auxilliary results_Malliavin differentiability} and similar arguments as in the proof of Proposition~\ref{Existence result for the FBSDE}, we obtain that the BSDE \eqref{E9_Proposition 3.9} admits a unique solution $(D_{t}(Y_{t}),D_{t}(Z_{t}))_{t\in[0,T]}\in\mathcal{S}_{\mathbb{Q}_{0}}^{2}\times\mathbb{H}_{\mathbb{Q}_{0}}^{2}$, with the expectation representation of the first component $(D_{t}(Y_{t}))_{t\in[0,T]}$ being given by
	{\small
		\begin{align}\label{E10_Proposition 3.9}
			D_{t}(Y_{t})&=\mathbb{E}^{\mathbb{Q}_{0}}\Big[D_{t}(-K_{T}e^{-\int_{0}^{T}r_{s}\mathrm{d}s})+\int_{t}^{T}\Big(\widetilde{x}D_{t}(\alpha_{s}\varphi_{s})\mathrm{d}s-Z_{s}^{\top}D_{t}(\eta_{s})\Big)\mathrm{d}s~|\mathcal{F}_{t}^{r,S}\Big].
	\end{align}}
	Using successively \eqref{E4_Proposition 3.9}, the fact that $\tilde{Z}_{t}=Z_{t},~t\in[0,T]$, and It{\^o} isometry we have
	\begin{align}\label{E11_Proposition 3.9}
		&\mathbb{E}^{\mathbb{Q}_{0}}\Big[\big(-K_{T}e^{-\int_{0}^{T}r_{s}\mathrm{d}s}+\widetilde{x}\int_{0}^{T}\alpha_{s}\varphi_{s}\mathrm{d}s\big)\int_{t}^{T}D_{t}(\eta_{s})\mathrm{d}I_{s}^{\mathbb{Q}_{0}}~|\mathcal{F}_{t}^{r,S}\Big]\notag\\
		&=\mathbb{E}^{\mathbb{Q}_{0}}\Big[\int_{t}^{T}Z_{s}^{\top}D_{t}(\eta_{s})\mathrm{d}s~|\mathcal{F}_{t}^{r,S}\Big].
	\end{align}
	Substituting \eqref{E11_Proposition 3.9} into \eqref{E10_Proposition 3.9} and using the linearity of the operator $D_{t}(\cdot)$ we obtain
	\begin{align}\label{E12_Proposition 3.9}
		D_{t}(Y_{t})&=\mathbb{E}^{\mathbb{Q}_{0}}\Big[D_{t}\big(-K_{T}e^{-\int_{0}^{T}r_{s}\mathrm{d}s}+\widetilde{x}\int_{0}^{T}\alpha_{s}\varphi_{s}\mathrm{d}s\big)\notag\\
		&\phantom{XXXXx}-\big(-K_{T}e^{-\int_{0}^{T}r_{s}\mathrm{d}s}+\widetilde{x}\int_{0}^{T}\alpha_{s}\varphi_{s}\mathrm{d}s\big)\int_{t}^{T}D_{t}(\eta_{s})\mathrm{d}I_{s}^{\mathbb{Q}_{0}}~|\mathcal{F}_{t}^{r,S}\Big].
	\end{align}
	Hence, comparing \eqref{E8_Proposition 3.9} and \eqref{E12_Proposition 3.9}, we deduce that $Z_{t}=D_{t}(Y_{t})$ for $t\in[0,T]$.
\end{proof}

\section{Utility loss}\label{Welfare loss}
In this section, we quantify the utility loss incurred when investors ignore the information revealed through learning about the market price of risk $R$. More precisely, instead of updating their estimate of $R$ through the filtering process $(m_{t})_{t\in[0,T]}$, investors behave as if the market price of risk were constantly equal to its long-run mean $\mu_{R}$. Following \cite{escobar2016portfolio}, we measure this utility loss as a percentage of the investor's initial wealth. Specifically, for $L\in(0,1)$, we define the utility loss by the indifference relation $\mathcal{V}(x(1-L))=\mathcal{V}(x)$, where $\mathcal{V}(x(1-L))$ denotes the value function associated with the optimisation problem \eqref{Main problem} when the initial wealth is reduced to $x(1-L)$, and $\mathcal{V}(x)$ denotes the value function corresponding to the optimisation problem in which the filtered estimate is replaced by the constant value $m_{t}=\mu_{R},\,t\in[0,T]$. Using Theorem \ref{Mainresult1}, we obtain
\begin{align*}
	L&=1-\frac{1}{x}\Bigg(\frac{1-\mathbb{E}\big[\int_{0}^{T}H_{s}\alpha_{s}\varphi_{s}\mathrm{d}s~|\mathcal{F}_{0}^{r,S}\big]}{1-\mathbb{E}\big[\int_{0}^{T}H_{s}^{0}\alpha_{s}^{0}\varphi_{s}^{0}\mathrm{d}s~|\mathcal{F}_{0}^{r,S}\big]}\Big(x-\mathbb{E}\big[H_{T}^{0}K_{T}e^{-\int_{0}^{T}r_{s}\mathrm{d}s}~|\mathcal{F}_{0}^{r,S}\big]\Big)\notag\\
	&\phantom{XXXxxx}+\mathbb{E}\big[H_{T}K_{T}e^{-\int_{0}^{T}r_{s}\mathrm{d}s}~|\mathcal{F}_{0}^{r,S}\big]\Bigg),
\end{align*}
where $H^{0},\alpha^{0}$, and $\varphi^{0}$ are defined by \eqref{Auxiliary processes} for $m_{t}=\mu_{R},~t\in[0,T]$.

In the remainder of this chapter, we restrict our attention to the case of a non-negative constant liability $K$, that is, $K_{T}=K$. Under this specification, it is straightforward to verify that Assumptions \ref{Change of measure}$.(ii)$ and \ref{Conditions to be Malliavin differentiable} are automatically satisfied.

\begin{figure}[h!]
	\centering
	\subfigure{
		\includegraphics[width=0.45\textwidth]{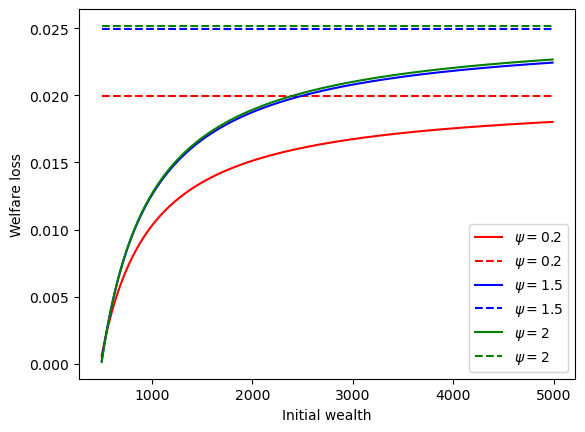}
		\label{fig:Effect of psi on the welfare loss}
	}
	\subfigure{
		\includegraphics[width=0.45\textwidth]{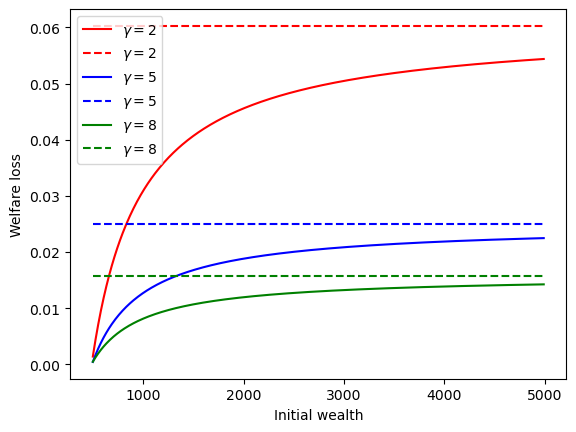}
		\label{fig:Effect of gamma on the welfare loss}
	}
	\caption{Welfare loss $L$. Both figures use $\delta=0.08$. The left panel uses $\gamma=5$, and the right panel takes $\psi=1.5$. The solid lines represent the cases where  $K=500$ and the dashed lines the cases where $K=0$.}
	\label{fig:Effect of gamma on the optimal strategies}
\end{figure}

To derive a sufficient condition ensuring that Assumption \ref{Change of measure}$.(i)$ holds, we introduce the quantities $\sigma_{m}^{2}(t):=\big(\sigma_{R}\rho_{RS}+\beta v_{t}\big)^{2}+\big(\sigma_{R}\rho_{Rr} -\rho_{rS}\beta(1-\rho_{rS}^{2})^{-\frac{1}{2}}v_{t}\big)^{2},~\Delta(t):=2\sigma_{m}^{2}(t)\zeta-\kappa_{R}^{2}$, $b_{max}:=\max_{t\in[0,T]}\sigma_{m}^{2}(t)$ and $\Delta_{max}:=2b_{max}\zeta-\kappa_{R}^{2}$, with $\zeta:=2q_{max}\beta^{2}(1-\rho_{rS}^{2})^{-1}$. Moreover, in the following results, we further specialise to the case $q=2$ in Assumption \ref{Change of measure}. Finally, we denote by \textsl{pi} the mathematical constant such that $180^{\circ}=\textsl{pi}$ radians.
\begin{prop}\label{Sufficient conditions for the assumptions to hold}
	Assume $\Delta_{max}>0,~T<\left(\textsl{pi}-\arctan(\sqrt{\Delta_{max}}/\kappa_{R})\right)/\sqrt{\Delta_{max}}$ or $\Delta_{max}\le0$ hold. Then, Assumption \ref{Change of measure}$.(i)$ is satisfied.
\end{prop}
\begin{proof}
	We define $u(t,m):=\mathbb{E}\Big[\exp\Big(\frac{2q_{max}\beta^{2}}{1-\rho_{rS}^{2}}\int_{t}^{T}m_{s}^{2}\mathrm{d}s\Big)~\big|m_{t}=m\Big]$. Then, using an exponential quadratic ansatz on $u$, the proof follows from Feynman-K{\v a}c formula and a comparison theorem for Riccati equations. See Appendix \ref{Appendix D} for a complete proof.
\end{proof}
In the numerical illustrations, we consider the following parameter values: $\kappa_{r}=0.5,\kappa_{R}=1.5,\mu_{r}=0.02,\mu_{R}=\phi_{B}=\rho_{rS}=0,\sigma_{r}=-0.03,\sigma_{R}=0.2,\beta=4,\rho_{RS}=-0.95,\rho_{Rr}=0.1$ and $T=1$. (All comparative statistics are produced using a Monte Carlo simulation of $1000000$ paths and averaging them). Figure \ref{fig:Effect of gamma on the optimal strategies} illustrates that, in the presence of liabilities, the welfare loss is increasing in the investors' initial wealth, whereas it is independent of initial wealth when liabilities are absent. Furthermore, the welfare loss decreases with the coefficient of relative risk aversion $\gamma$ and increases with the elasticity of intertemporal substitution $\psi$, indicating that the value of learning is greater for investors who are more willing to substitute consumption across time. This observation also explains why the welfare loss is smaller under CRRA utility than under Epstein--Zin utility: under the CRRA restriction $\gamma\psi=1$, a higher degree of risk aversion is necessarily associated with a lower elasticity of intertemporal substitution, thereby reducing the value of learning.

\appendix
\section{Integrability of $\alpha\varphi$}\label{Appendix A}
Recall that expectations are taken under the observation filtration $\mathbb{F}^{r,S}$. For $\gamma,\psi>1$, using Jensen's inequality we have
\begin{align*}
	&\mathbb{E}^{\mathbb{Q}_{0}}\Big[\Big(\int_{0}^{T}|\alpha_{s}\varphi_{s}|^{2}\mathrm{d}s\Big)^{\frac{2q}{2}}\Big]\le T^{q-1}\mathbb{E}^{\mathbb{Q}_{0}}\Big[\int_{0}^{T} |\alpha_{s}\varphi_{s}|^{2q}\mathrm{d}s\Big]\notag\\
	&=T^{q-1}\mathbb{E}\Big[\int_{0}^{T}e^{-2q\int_{0}^{s}r_{u}\mathrm{d}u}\Big(\frac{\delta^{\psi}}{\psi-1}+r_{s}+\frac{1}{2\gamma}\|\eta_{s}\|^{2}{-\frac{\delta\theta}{1-\gamma}}\Big)^{2q}\notag\\
	&\phantom{x}\times\exp\Big(\int_{0}^{s}\Big(-2q\frac{\delta^{\psi}\psi}{\psi-1}+\frac{q}{\gamma}\|\eta_{u}\|^{2}{+2q\frac{\delta\theta}{1-\gamma}}\Big)\mathrm{d}u\Big) \notag\\
	&\phantom{x}\times\exp\Big(-\frac{1}{2}\int_{0}^{s} \|\eta_{u}\|^{2}\mathrm{d}u-\int_{0}^{s}\eta_{u}^{\top}\mathrm{d}I_{u}\Big)\exp\Big(-\frac{q}{\gamma^{2}}\int_{0}^{s}\|\eta_{u}\|^{2}\mathrm{d}u+\frac{2q}{\gamma}\int_{0}^{s}\eta_{u}^{\top}\mathrm{d}I_{u}\Big)\mathrm{d}s\Big]\notag\\
	&=T^{q-1}\mathbb{E}\Big[\int_{0}^{T}e^{-2q\int_{0}^{s}r_{u}\mathrm{d}u}\Big(\frac{\delta^{\psi}}{\psi-1}+r_{s}+\frac{1}{2\gamma}\|\eta_{s}\|^{2}{-\frac{\delta\theta}{1-\gamma}}\Big)^{2q}\notag\\
	&\phantom{x}\times\exp\Big(\int_{0}^{s}\Big(-2q\frac{\delta^{\psi}\psi}{\psi-1}+\frac{\gamma^{2}-5q\gamma+6q^{2}-q}{\gamma^{2}}\|\eta_{u}\|^{2}{+2q\frac{\delta\theta}{1-\gamma}}\Big)\mathrm{d}u\Big)\mathrm{d}s\Big]\notag\\
	&\phantom{x}\times\exp\Big(-\frac{1}{6}\Big(\frac{3(2q-\gamma)}{\gamma}\Big)^{2}\int_{0}^{s}\|\eta_{u}\|^{2}\mathrm{d}u+\frac{2q-\gamma}{\gamma}\int_{0}^{s}\eta_{u}^{\top}\mathrm{d}I_{u}\Big).
\end{align*}
Hence
\begin{align*}
	&\mathbb{E}^{\mathbb{Q}_{0}}\Big[\Big(\int_{0}^{T}|\alpha_{s}\varphi_{s}|^{2}\mathrm{d}s\Big)^{\frac{2q}{2}}\Big]\notag\\
	&\le\Big(\mathbb{E}\Big[\int_{0}^{T}e^{-6q\int_{0}^{s}r_{u}\mathrm{d}u}\Big(\frac{\delta^{\psi}}{\psi-1}+r_{s}+\frac{1}{2\gamma}\|\eta_{s}\|^{2}{-\frac{\delta\theta}{1-\gamma}}\Big)^{6q}\mathrm{d}s\Big]\Big)^{\frac{1}{3}}\notag\\
	&\phantom{x}\times\Big(\mathbb{E}\Big[\int_{0}^{T}\mathcal{E}\Big(\int\frac{3(2q-\gamma)}{\gamma}\eta^{\top}\mathrm{d}I\Big)_{s}\mathrm{d}s\Big]\Big)^{\frac{1}{3}}\notag\\
	&\phantom{x}\times\Big(\mathbb{E}\Big[\int_{0}^{T}\exp\Big(\int_{0}^{s}\Big(3\frac{\gamma^{2}-5q\gamma+6q^{2}-q}{\gamma^{2}}\|\eta_{u}\|^{2}\Big)\mathrm{d}u\Big)\mathrm{d}s\Big]\Big)^{\frac{1}{3}}\exp\Big({2q\frac{\delta\theta}{1-\gamma}}T\Big)T^{q-1}\notag\\
	&\le\Big(\mathbb{E}\Big[\int_{0}^{T}e^{-12q\int_{0}^{s}r_{u}\mathrm{d}u}\mathrm{d}s\Big]\Big)^{\frac{1}{6}}\Big(\mathbb{E}\Big[\int_{0}^{T}\Big(\frac{\delta^{\psi}}{\psi-1}+r_{s}+\frac{1}{2\gamma}\|\eta_{s}\|^{2}{-\frac{\delta\theta}{1-\gamma}}\Big)^{12q}\mathrm{d}s\Big]\Big)^{\frac{1}{6}}\notag\\
	&\phantom{x}\times\Big(\mathbb{E}\Big[\int_{0}^{T}\mathcal{E}\Big(\int\frac{3(2q-\gamma)}{\gamma}\eta^{\top}\mathrm{d}I\Big)_{s}\mathrm{d}s\Big]\Big)^{\frac{1}{3}}\exp\Big({2q\frac{\delta\theta}{1-\gamma}}T\Big)T^{q-1}\notag\\
	&\phantom{x}\times\Big(\mathbb{E}\Big[\int_{0}^{T}\exp\Big(\int_{0}^{s}\Big(3\frac{\gamma^{2}-5q\gamma+6q^{2}-q} {\gamma^{2}}\|\eta_{u}\|^{2}\Big)\mathrm{d}u\Big)\mathrm{d}s\Big]\Big)^{\frac{1}{3}}\notag\\
	&\le\Big(\mathbb{E}\Big[\int_{0}^{T}\Big(3^{12q-1}\Big(\frac{\delta^{\psi}}{\psi-1}-\frac{\delta\theta}{1-\gamma}\Big)^{12q}+3^{12q-1}|r_{s}|^{12q}+\frac{3^{12q-1}}{(2\gamma)^{12q}}\|\eta_{s}\|^{12q}\Big)\mathrm{d}s\Big]\Big)^{\frac{1}{6}}\notag\\
	&\phantom{x}\times\Big(\mathbb{E}\Big[\int_{0}^{T}e^{-12q\int_{0}^{s}r_{u}\mathrm{d}u}\mathrm{d}s\Big]\Big)^{\frac{1}{6}}\exp\Big({2q\frac{\delta\theta}{1-\gamma}}T\Big)T^{q-\frac{2}{3}}\notag\\
	&\phantom{x}\times\Big(\mathbb{E}\Big[\int_{0}^{T}\exp\Big(3\frac{\gamma^{2}-5q\gamma+6q^{2}-q}{\gamma^{2}}\int_{0}^{s}\|\eta_{u}\|^{2}\mathrm{d}u\Big)\mathrm{d}s\Big]\Big)^{\frac{1}{3}}\notag\\
	&<\infty,
\end{align*}
where the first inequality holds due to H{\"o}lder's inequality, the facts that $-2q\frac{\delta^{\psi}\psi}{\psi-1}<0$ and $\exp\Big({2q\frac{\delta\theta}{1-\gamma}}s\Big)\le\exp\Big({2q\frac{\delta\theta}{1-\gamma}}T\Big)$ for $s\in[0,T]$, the second inequality follows from Cauchy-Schwarz inequality, the third inequality comes from the convex inequality $(a+b+c)^{12q}\le3^{12q-1}\big(a^{12q}+b^{12q}+c^{12q}\big)$ and Assumption~\ref{Change of measure}.$(i)$, and the last inequality holds due to Assumption~\ref{Change of measure}.$(i)$ and the facts that $(r_{t})_{t\in[0,T]}$ is still an OU processs in the observation filtration $\mathbb{F}^{r,S}$ (see \eqref{Filtered interest rate}) and $\|\eta_{t}\|^{2}$ is the sum of power of the OU process $(m_{t})_{t\in[0,T]}$. Hence, $\alpha\varphi\in\mathbb{H}_{\mathbb{Q}_{0}}^{2q}$ for all $q>1$.

\section{Proof of Proposition \ref{Existence result for the FBSDE}}\label{Appendix B}
The idea of the proof is to apply a decoupling reduction method via which an explicit solution to the coupled linear FBSDE \eqref{Optimal FBSDE} with unbounded random coefficients can be obtained. The key difference from existing works such as \cite{ma1994solving,yong1999linear,yong2006linear}, where the auxiliary processes $P$ and $p$ therein are solutions of BSDEs, is that our approach introduces two auxiliary processes $R$ and $\widetilde{\varphi}$ defined via forward SDEs. We then express the forward component of the FBSDE in terms of the backward process $Y$ and the forward processes $R$ and $\widetilde{\varphi}$. More precisely, we assume
\begin{align}
	X_{t}^{c^{*},\pi^{*}}&=R_{t}Y_{t}+\widetilde{\varphi}_{t}~\text{ for all }t\in[0,T].
\end{align}

Owing to the structure of the optimal value process given in \eqref{Ansatz_Optimal process}, we set
\begin{align}
	R_{t}=-e^{\int_{0}^{t}r_{s}\mathrm{d}s}~\text{ for all }t\in[0,T].
\end{align} 
Hence, we are looking for a solution $(X^{c^{*},\pi^{*}},Y,Z)$ such that
\begin{align}\label{Ansatz on X}
	X_{t}^{c^{*},\pi^{*}}&=-e^{\int_{0}^{t}r_{s}\mathrm{d}s}Y_{t}+\widetilde{\varphi}_{t}~\text{ for all }t\in[0,T],
\end{align}
where $\widetilde{\varphi}:[0,T]\times\Omega\rightarrow\mathbb{R}$ is an adapted process solution to the SDE
\begin{align}\label{Dynamics of varphi}
	\mathrm{d}\widetilde{\varphi}_{t}&=\widetilde{\alpha}_{t}\mathrm{d}t+\widetilde{\beta}_{t}\mathrm{d}I_{t},~~\widetilde{\varphi}_{0}\in\mathbb{R},
\end{align}
with the constant $\widetilde{\varphi}_{0}$ and the processes $\widetilde{\alpha}$ and $\widetilde{\beta}$ to be determined. Applying It{\^o}'s formula to the process $-e^{\int_{0}^{t}r_{s}\mathrm{d}s}Y_{t}+\widetilde{\varphi}_{t},\,t\in[0,T]$, and using \eqref{Dynamics of varphi} and the BSDE in \eqref{Optimal FBSDE} we obtain
\begin{align}\label{WealthOptimal_1}
	&\mathrm{d}\big(-e^{\int_{0}^{t}r_{s}\mathrm{d}s}Y_{t}+\widetilde{\varphi}_{t}\big)\notag\\
	&=\Big[-r_{t}e^{\int_{0}^{t}r_{s}\mathrm{d}s}Y_{t}+\Big(\frac{\delta^{\psi}}{\psi-1}+r_{t}+\frac{1}{2\gamma}\|\eta_{t}\|^{2}{-\frac{\delta\theta}{1-\gamma}}\Big)(X_{t}^{c^{*},\pi^{*}}+e^{\int_{0}^{t}r_{s}\mathrm{d}s}Y_{t})\notag\\
	&\phantom{xx}-e^{\int_{0}^{t}r_{s}\mathrm{d}s}Z_{t}^{\top}\eta_{t}+\widetilde{\alpha}_{t}\Big]\mathrm{d}t +\Big(-e^{\int_{0}^{t}r_{s}\mathrm{d}s}Z_{t}^{\top}+\widetilde{\beta}_{t}\Big)\mathrm{d}I_{t}\notag\\
	&=\Big[-e^{\int_{0}^{t}r_{s}\mathrm{d}s}Z_{t}^{\top} \eta_{t}+r_{t}X_{t}^{c^{*},\pi^{*}}+\big(-\delta^{\psi}+\frac{1}{\gamma}\|\eta_{t}\|^{2}\big)(X_{t}^{c^{*},\pi^{*}}+e^{\int_{0}^{t}r_{s}\mathrm{d}s}Y_{t})\notag\\
	&\phantom{xx}+\Big(\frac{\delta^{\psi}\psi}{\psi-1}-\frac{1}{2\gamma}\|\eta_{t}\|^{2}{-\frac{\delta\theta}{1-\gamma}}\Big)(X_{t}^{c^{*},\pi^{*}}+e^{\int_{0}^{t}r_{s}\mathrm{d}s}Y_{t})+\widetilde{\alpha}_{t}\Big]\mathrm{d}t\notag\\
	&\phantom{xx}+\Big[-e^{\int_{0}^{t}r_{s}\mathrm{d}s}Z_{t}^{\top}+\widetilde{\beta}_{t}+\frac{1}{\gamma}(X_{t}^{c^{*},\pi^{*}}+e^{\int_{0}^{t}r_{s}\mathrm{d}s}Y_{t})\eta_{t}^{\top}-\frac{1}{\gamma}(X_{t}^{c^{*},\pi^{*}}+e^{\int_{0}^{t}r_{s}\mathrm{d}s}Y_{t})\eta_{t}^{\top}\Big]\mathrm{d}I_{t}.
\end{align}
Plugging the expression of $\widetilde{\varphi}$ (see \eqref{Ansatz on X}) in the second line of the second equation in \eqref{WealthOptimal_1}, we deduce from \eqref{Ansatz on X} that
\begin{align}\label{WealthOptimal_2}
	\mathrm{d}X_{t}^{c^{*},\pi^{*}}&=\Big[r_{t}X_{t}^{c^{*},\pi^{*}}+\big(-\delta^{\psi}+\frac{1}{\gamma}\|\eta_{t}\|^{2}\big)(X_{t}^{c^{*},\pi^{*}}+e^{\int_{0}^{t}r_{s}\mathrm{d}s}Y_{t})-e^{\int_{0}^{t}r_{s}\mathrm{d}s}Z_{t}^{\top} \eta_{t}\notag\\
	&\phantom{xx}+\Big(\frac{\delta^{\psi}\psi}{\psi-1}-\frac{1}{2\gamma}\|\eta_{t}\|^{2}{-\frac{\delta\theta}{1-\gamma}}\Big)\widetilde{\varphi}_{t}+\widetilde{\alpha}_{t}\Big]\mathrm{d}t\notag\\
	&\phantom{xx}+\Big[\frac{1}{\gamma}(X_{t}^{c^{*},\pi^{*}}+e^{\int_{0}^{t}r_{s}\mathrm{d}s}Y_{t})\eta_{t}^{\top}-e^{\int_{0}^{t}r_{s}\mathrm{d}s}Z_{t}^{\top}-\frac{1}{\gamma}\widetilde{\varphi}_{t}\eta_{t}^{\top}+\widetilde{\beta}_{t}\Big]\mathrm{d}I_{t}.
\end{align}
Hence, comparing \eqref{WealthOptimal_2} to the Forward SDE in \eqref{Optimal FBSDE}, we have
\begin{align}
	\widetilde{\alpha}_{t}&=\Big(-\frac{\delta^{\psi}\psi}{\psi-1}+\frac{1}{2\gamma}\|\eta_{t}\|^{2}{+\frac{\delta\theta}{1-\gamma}}\Big)\widetilde{\varphi}_{t}~\text{ and }~
	\widetilde{\beta}_{t}=\frac{1}{\gamma}\widetilde{\varphi}_{t}\eta_{t}^{\top}~\text{ for all } ~t\in[0,T].
\end{align}
Thus, $\widetilde{\varphi}$ solves the SDE
\begin{align}\label{Dynamics of the auxiliary process_1}
	\mathrm{d}\widetilde{\varphi}_{t}&=\widetilde{\varphi}_{t}\Big(-\frac{\delta^{\psi}\psi}{\psi-1}+\frac{1}{2\gamma}\|\eta_{t}\|^{2}{+\frac{\delta\theta}{1-\gamma}}\Big)\mathrm{d}t+\frac{1}{\gamma}\widetilde{\varphi}_{t} \eta_{t}^{\top}\mathrm{d}I_{t},\text{ with }~\widetilde{\varphi}_{0}\in\mathbb{R}.
\end{align}
The following result gives sufficient conditions for the \textit{solvability} of the FBSDE \eqref{Optimal FBSDE}.
\begin{lemm}\label{Solvability of the FBSDE}
	Let Assumption \ref{Change of measure} hold. Let $\alpha$ be the process given in \eqref{Auxiliary processes}. Assume that the constant $\widetilde{x}$ is finite and the SDE
	\begin{align}\label{Dynamics of the auxiliary process_1.0}
		\mathrm{d}\widetilde{\varphi}_{t}&=\widetilde{\varphi}_{t}\Big(-\frac{\delta^{\psi}\psi}{\psi-1}+\frac{1}{2\gamma}\|\eta_{t}\|^{2}{+\frac{\delta\theta}{1-\gamma}}\Big)\mathrm{d}t+\frac{1}{\gamma}\widetilde{\varphi}_{t}\eta_{t}^{\top}\mathrm{d}I_{t},\text{ with }~\widetilde{\varphi}_{0}=\widetilde{x},
	\end{align}
	admits a unique strong solution $\widetilde{\varphi}\in\mathcal{H}_{\mathbb{P}}^{q},~q>1$, and the BSDE
	\begin{align}\label{Auxiliary BSDE_2.0}
		\mathrm{d}Y_{t}&=-\Big(\alpha_{t}\widetilde{\varphi}_{t} -Z_{t}^{\top}\eta_{t}\Big)\mathrm{d}t+Z_{t}^{\top}\mathrm{d}I_{t},~~Y_{T}=-e^{-\int_{0}^{T}r_{s}\mathrm{d}s}K_{T}
	\end{align}
	admits a unique solution $(Y,Z)\in\mathcal{S}_{\mathbb{P}}^{q}\times\mathbb{H}_{\mathbb{P}}^{q},~q>1$. Then the FBSDE \eqref{Optimal FBSDE} admits a unique solution $(X^{c^{*},\pi^{*}},Y,Z)\in\mathcal{H}_{\mathbb{P}}^{q'}\times\mathcal{S}_{\mathbb{P}}^{q}\times\mathbb{H}_{\mathbb{P}}^{q},~q>q'>1$, such that the representation \eqref{Ansatz on X} holds.
\end{lemm}
\begin{proof}
	Substituting in \eqref{Auxiliary BSDE_2.0} the adapted processes $\alpha$ and $\widetilde{\varphi}$ by their expression given by \eqref{Auxiliary processes} and \eqref{Ansatz on X}, respectively, we obtain the BSDE in \eqref{Optimal FBSDE}. Then, applying It{\^o}'s formula to \eqref{Ansatz on X} we obtain the Forward SDE in \eqref{Optimal FBSDE}. Moreover, for all $q'\in(1,q)$, using the convex inequality $(a+b)^{q'}\le a^{q'}+b^{q'}$, Young inequality, Assumption \ref{Change of measure}$.(i)$ and the fact that $\widetilde{\varphi}\in\mathcal{H}_{\mathbb{P}}^{q},\,q>q'$, we have
	\begin{align*}
		&\mathbb{E}\Big[\int_{0}^{T}|X_{t}^{c^{*},\pi^{*}}|^{q'}\mathrm{d}t~\big|\mathcal{F}_{0}^{r,S}\Big]\notag\\
		&\le\mathbb{E}\Big[\sup_{t\le T} |Y_{t}|^{q}\Big]+\mathbb{E}\Big[\int_{0}^{T}e^{-\frac{qq'}{q-q'}\int_{0}^{t}r_{s}\mathrm{d}s}\mathrm{d}t~\big|\mathcal{F}_{0}^{r,S}\Big]+\mathbb{E}\Big[\int_{0}^{T}|\widetilde{\varphi}_{t}|^{q'}\mathrm{d}t~\big|\mathcal{F}_{0}^{r,S}\Big]\notag\\
		&<\infty.
	\end{align*}
	Hence, $(X^{c^{*},\pi^{*}},Y,Z)\in\mathcal{H}_{\mathbb{P}}^{q'}\times\mathcal{S}_{\mathbb{P}}^{q}\times\mathbb{H}_{\mathbb{P}}^{q},~q>1$, is an adapted solution to the FBSDE \eqref{Optimal FBSDE}. The uniquenss follows from the uniqueness of $\widetilde{\varphi}$ (the solution of the SDE \eqref{Dynamics of the auxiliary process_1.0}) and the uniqueness of the pair $(Y,Z)$ (the solution of the BSDE \eqref{Auxiliary BSDE_2.0}).
\end{proof}

Thanks to Lemma \ref{Solvability of the FBSDE}, to establish the existence of a (unique) solution to the FBSDE system \eqref{Optimal FBSDE}, it suffices to show that the decoupled FBSDE \eqref{Dynamics of the auxiliary process_1.0}–\eqref{Auxiliary BSDE_2.0} admits a unique solution. We start by showing that under Assumption \ref{Change of measure} the SDE has a unique solution.
\begin{lemm}\label{Solution to the auxiliary SDE}
	Let Assumption \ref{Change of measure} hold. Then, the constant $\widetilde{x}$ is finite and the SDE \eqref{Dynamics of the auxiliary process_1.0} admits a unique strong solution $\widetilde{\varphi}\in\mathcal{H}_{\mathbb{P}}^{q},~q>1$, given by
	\begin{align}\label{Auxiliary process_1}
		\widetilde{\varphi}_{t}&=\widetilde{x}\exp\Big(\int_{0}^{t}\Big(-\frac{\delta^{\psi}\psi}{\psi-1}+\frac{\gamma-1}{2\gamma^{2}}\|\eta_{s}\|^{2}{+\frac{\delta\theta}{1-\gamma}}\Big)\mathrm{d}s+\int_{0}^{t}\frac{1}{\gamma}\eta_{s}^{\top}\mathrm{d}I_{s}\Big),~0\le t\le T.
	\end{align}
\end{lemm}
\begin{proof}
	To show that $\widetilde{\varphi}$ from \eqref{Auxiliary process_1} solves \eqref{Dynamics of the auxiliary process_1.0} it suffices to apply It{\^o}'s formula on $\widetilde{\varphi}$. Now, suppose there exists another solution $(\overline{\varphi}_{t})_{t\in[0,T]}$ to \eqref{Dynamics of the auxiliary process_1.0}. We define the processes
	\begin{align*}
		\widetilde{\beta}_{t}&:=\exp\Big(\int_{0}^{t}\Big(\frac{\delta^{\psi}\psi}{\psi-1}-\frac{\gamma+1}{2\gamma^{2}}\|\eta_{s}\|^{2}{-\frac{\delta\theta}{1-\gamma}}\Big)\mathrm{d}s -\int_{0}^{t}\frac{1}{\gamma}\eta_{s}^{\top}\mathrm{d}I_{s}\Big)\notag\\
		\widetilde{\Delta}_{t}&:=\widetilde{\beta}_{t}\left(\widetilde{\varphi}_{t}-\overline{\varphi}_{t}\right),~0\le t\le T.
	\end{align*}
	Applying It{\^o}'s formula on $\widetilde{\Delta}$ we obtain the uniqueness of the solution to the SDE~\eqref{Dynamics of the auxiliary process_1.0}.
	
	Finally, proving that $\widetilde{\varphi}\in\mathcal{H}_{\mathbb{P}}^{q},~q>1$, follows directly from Assumption~\ref{Change of measure}. Indeed, using the H{\"o}lder's inequality, Assumption~\ref{Change of measure}$.(i)$,  with $\frac{4q^{2}}{\gamma^{2}}<\frac{9(2q-\gamma)^{2}}{2\gamma^{2}}\le q_{max}$ and $\frac{2q^{2}+q(\gamma-1)}{\gamma^{2}}\le q_{max}$, we have
	\begin{align*}
		&\mathbb{E}\Big[\int_{0}^{T}|\widetilde{\varphi}_{t}|^{q}\mathrm{d}t~\big|\mathcal{F}_{0}^{r,S}\Big]\notag\\
		&=\mathbb{E}\Big[\int_{0}^{T}\exp\Big(\int_{0}^{t}\Big(-q\frac{\delta^{\psi}\psi}{\psi-1}+\frac{2q^{2}+q(\gamma-1)}{2\gamma^{2}}\|\eta_{s}\|^{2}{+q\frac{\delta\theta}{1-\gamma}}\Big)\mathrm{d}s\Big)\notag\\
		&\times\exp\Big(-\frac{1}{4}\Big(\frac{2q}{\gamma}\Big)^{2}\int_{0}^{t}\|\eta_{s}\|^{2}\mathrm{d}s+\frac{q}{\gamma}\int_{0}^{t}\eta_{s}^{\top}\mathrm{d}I_{s}\Big)\mathrm{d}t~\big|\mathcal{F}_{0}^{r,S}\Big]\notag\\
		&\le\Big(\mathbb{E}\Big[\int_{0}^{T}\exp\Big(\int_{0}^{t}\Big(-2q\frac{\delta^{\psi}\psi}{\psi-1}+\frac{2q^{2}+q(\gamma-1)}{\gamma^{2}}\|\eta_{s}\|^{2}{+2q\frac{\delta\theta}{1-\gamma}}\Big)\mathrm{d}s\Big)~\big|\mathcal{F}_{0}^{r,S}\Big]\Big)^{\frac{1}{2}}\notag\\
		&\times\Big(\mathbb{E}\Big[\int_{0}^{T}\mathcal{E}\Big(\int\frac{2q}{\gamma}\eta^{\top}\mathrm{d}W\Big)_{t}\mathrm{d}t~\big|\mathcal{F}_{0}^{r,S}\Big]\Big)^{\frac{1}{2}}\notag\\
		&\le T^{\frac{3}{2}}\exp\Big(\Big|-2q\frac{\delta^{\psi}\psi}{\psi-1}+2q\frac{\delta\theta}{1-\gamma}\Big|T\Big)\Big(\mathbb{E}\Big[\exp\Big(\frac{2q^{2}+q(\gamma-1)}{\gamma^{2}}\int_{0}^{T}\|\eta_{s}\|^{2}\mathrm{d}s\Big)~\big|\mathcal{F}_{0}^{r,S}\Big]\Big)^{\frac{1}{2}}\notag\\
		&<\infty.
	\end{align*}
	
	Finally, we prove that the constant $\widetilde{x}$ is finite. By Assumption~\ref{Change of measure}$.(i)$, it suffices to show that $1-\mathbb{E}\big[\int_{0}^{T}H_{s}\alpha_{s}\varphi_{s}\mathrm{d}s\big]\ne0$. Indeed, recalling the expressions of $\alpha$ and $\varphi$ from \eqref{Auxiliary processes} , we have
	\begin{align}\label{E1_Proposition 3.2}
		&1-\mathbb{E}\Big[\int_{0}^{T}H_{s}\alpha_{s}\varphi_{s}\mathrm{d}s~\big|\mathcal{F}_{0}^{r,S}\Big]\notag\\
		&=\mathbb{E}^{\mathbb{Q}_{1}}\Big[\int_{0}^{T}\delta^{\psi}\exp\Big(\int_{0}^{s}\big(-\frac{\delta^{\psi}\psi}{\psi-1}-r_{u}-\frac{1}{2\gamma}\|\eta_{u}\|^{2}{+\frac{\delta\theta}{1-\gamma}}\big)\mathrm{d}u\Big)\mathrm{d}s~\big|\mathcal{F}_{0}^{r,S}\Big]\notag\\
		&\phantom{xx}+\mathbb{E}^{\mathbb{Q}_{1}}\Big[\exp\Big(\int_{0}^{T}\big(-\frac{\delta^{\psi}\psi}{\psi-1}-r_{u}-\frac{1}{2\gamma}\|\eta_{u}\|^{2}{+\frac{\delta\theta}{1-\gamma}}\big)\mathrm{d}u\Big)~\big|\mathcal{F}_{0}^{r,S}\Big]>0,
	\end{align}
	where $\mathbb{Q}_{1}$ is the probability measure equivalent to $\mathbb{P}$ and defined by\\ $\frac{\mathrm{d}\mathbb{Q}_{1}}{\mathrm{d}\mathbb{P}}\big|_{\mathcal{F}_{T}^{r,S}}:=\exp\left(-\frac{(1-\gamma)^{2}}{2\gamma^{2}}\int_{0}^{T}\|\eta_{s}\|^{2}\mathrm{d}s+\frac{1-\gamma}{\gamma}\int_{0}^{T}\eta_{s}^{\top}\mathrm{d}I_{s}\right)$. Note that $\mathbb{Q}_{1}$ is well-defined thanks to Assumption \ref{Change of measure}$.(i)$.
\end{proof}

The following result concerns the \textit{solvability} of the BSDE \eqref{Auxiliary BSDE_2.0}.
\begin{lemm}\label{Solution to the auxiliary BSDE}
	Let Assumption \ref{Change of measure} hold. Let $\alpha$ be the process given in \eqref{Auxiliary processes}. Let $\widetilde{\varphi}$ be the unique solution to \eqref{Dynamics of the auxiliary process_1.0}. Then, the BSDE \eqref{Auxiliary BSDE_2.0}	admits a unique solution $(Y,Z)\in\mathcal{S}_{\mathbb{P}}^{q}\times\mathbb{H}_{\mathbb{P}}^{q},~q>1$, and the expectation representation of the first component $Y$ is given by
	\begin{align}
		Y_{t}&=H_{t}^{-1}\mathbb{E}\Big[-H_{T}K_{T}e^{-\int_{0}^{T}r_{s}\mathrm{d}s}+\int_{t}^{T}H_{s}\alpha_{s}\widetilde{\varphi}_{s}\mathrm{d}s~\big|\mathcal{F}_{t}^{r,S}\Big],\,\,0\le t\le T,
	\end{align}
	with the process $H$ given in \eqref{Auxiliary processes}.
\end{lemm}
\begin{proof}
	Under the new measure $\mathbb{Q}_{0}$ (see Assumption \ref{Change of measure}$.(ii)$), we consider a pair $(\tilde{Y},\tilde{Z})$ satisfying the BSDE
	\begin{align}
		\mathrm{d}\tilde{Y}_{t}&=-\alpha_{t}\widetilde{\varphi}_{t}\mathrm{d}t+\tilde{Z}_{t}^{\top}\mathrm{d}I_{t}^{\mathbb{Q}_{0}}=-\Big(\alpha_{t}\widetilde{\varphi}_{t}-\tilde{Z}_{t}^{\top}\eta_{t}\Big)\mathrm{d}t+\tilde{Z}_{t}^{\top}\mathrm{d}I_{t},\label{E1_Proposition 3.3}
	\end{align}
	with $\tilde{Y}_{T}=-K_{T}e^{-\int_{0}^{T}r_{s}\mathrm{d}s}$. Using Appendix \ref{Appendix A} and \cite[Thm.~5.1]{el1997backward}, the BSDE \eqref{E1_Proposition 3.3} admits a unique solution $(\tilde{Y},\tilde{Z})\in{\mathcal{S}_{\mathbb{Q}_{0}}^{2q}}\times\mathbb{H}_{\mathbb{Q}_{0}}^{2q},~{q>1}$, with the expectation representation of the first component $\tilde{Y}$ being given by
	\begin{align}\label{E2_Proposition 3.3}
		\tilde{Y}_{t}&=\mathbb{E}^{\mathbb{Q}_{0}}\Big[-K_{T}e^{-\int_{0}^{T}r_{s}\mathrm{d}s}+\widetilde{x}\int_{t}^{T}\alpha_{s}\widetilde{\varphi}_{s}\mathrm{d}s~\big|\mathcal{F}_{t}^{r,S}\Big]\notag\\
		&=H_{t}^{-1}\mathbb{E}\Big[-H_{T}K_{T}e^{-\int_{0}^{T}r_{s}\mathrm{d}s}+\widetilde{x}\int_{t}^{T}H_{s}\alpha_{s}\widetilde{\varphi}_{s}\mathrm{d}s~\big|\mathcal{F}_{t}^{r,S}\Big],~0\le t\le T.
	\end{align}
	From \eqref{E1_Proposition 3.3} we deduce that $(\tilde{Y},\tilde{Z})$ is also the unique solution to \eqref{Auxiliary BSDE_2.0}. Moreover, using  repeatedly Cauchy-Schwarz inequality we obtain
	\begin{align*}
		\mathbb{E}\Big[\sup_{t\le T}|\tilde{Y}_{t}|^{q}\mathrm{d}t~\big|\mathcal{F}_{0}^{r,S}\Big]
		&\le\Big(\mathbb{E}^{\mathbb{Q}_{0}}\Big[H_{T}^{-2}~\big|\mathcal{F}_{0}^{r,S}\Big]\Big)^{\frac{1}{2}}\Big(\mathbb{E}^{\mathbb{Q}_{0}}\Big[\sup_{t\le T}|\tilde{Y}_{t}|^{2q}\mathrm{d}t~\big|\mathcal{F}_{0}^{r,S}\Big]\Big)^{\frac{1}{2}}\notag\\
		&\le\Big(\mathbb{E}\Big[\exp\Big(3\int_{0}^{T}\|\eta_{s}\|^{2}\mathrm{d}s\Big)~\big|\mathcal{F}_{0}^{r,S}\Big]\Big)^{\frac{1}{4}}\Big(\mathbb{E}^{\mathbb{Q}_{0}}\Big[\sup_{t\le T}|\tilde{Y}_{t}|^{2q}\mathrm{d}t~\big|\mathcal{F}_{0}^{r,S}\Big]\Big)^{\frac{1}{2}}\notag\\
		&<\infty,
	\end{align*}
	where the last inequality holds due to Assumption~\ref{Change of measure}$.(i)$ and the fact that $\tilde{Y}\in{\mathcal{S}_{\mathbb{Q}_{0}}^{2q}}$. Using similar arguments and the fact that $\tilde{Z}\in\mathbb{H}_{\mathbb{Q}_{0}}^{2q}$, we have
	\begin{align*}
		&\mathbb{E}\Big[\Big(\int_{0}^{T}|\tilde{Z}_{s}|^{2}\mathrm{d}s\Big)^{\frac{q}{2}}~\big|\mathcal{F}_{0}^{r,S}\Big]\notag\\
		&\le\Big(\mathbb{E}\Big[\exp\Big(3\int_{0}^{T}\|\eta_{s}\|^{2}\mathrm{d}s\Big)~\big|\mathcal{F}_{0}^{r,S}\Big]\Big)^{\frac{1}{4}}\Big(\mathbb{E}^{\mathbb{Q}_{0}}\Big[\Big(\int_{0}^{T}|\tilde{Z}_{s} |^{2}\mathrm{d}s\Big)^{q}~\big|\mathcal{F}_{0}^{r,S}\Big]\Big)^{\frac{1}{2}}\notag\\
		&<\infty.
	\end{align*}
\end{proof}

We can now prove that the linear FBSDE \eqref{Optimal FBSDE} admits a solution.
\begin{proof}[Proof of Proposition \ref{Existence result for the FBSDE}]
	From \eqref{Auxiliary processes} and \eqref{Auxiliary process_1}, we deduce that $\widetilde{\varphi}_{t}=\widetilde{x}\varphi_{t},\,t\in[0,T]$. The remaining of the proof follows directly from Lemmas \ref{Solvability of the FBSDE}, \ref{Solution to the auxiliary SDE} and \ref{Solution to the auxiliary BSDE}.
\end{proof}

\section{Proof of Proposition \ref{Auxilliary results_Malliavin differentiability}}\label{Appendix C}
Recall that expectations are taken under the observation filtration $\mathbb{F}^{r,S}$. We first state and prove an intermediate result (Lemma \ref{Solution of the variance equation}). Lemma \ref{Solution of the variance equation} gives the expression of the solution of the Riccati equation given in Proposition \ref{Conditional mean-variance pair} and presents the bounds of such solution.
\begin{lemm}\label{Solution of the variance equation}
	For $\beta\ne0$, the solution $v$ to the Riccati equation
	\begin{align}\label{Variance equation}
		\mathrm{d}v_{t}&=\Big(\sigma_{R}^{2}-2\kappa_{R}v_{t}-\big(\sigma_{R}\rho_{RS}+\beta v_{t}\big)^{2}-\big(\sigma_{R}\rho_{Rr}-\rho_{rS}\beta (1-\rho_{rS}^{2})^{-\frac{1}{2}}v_{t}\big)^{2}\Big)\mathrm{d}t,~v_{0}\ge0
	\end{align}
	is given by
	\begin{align}\label{Solution of the variance equation_Expression}
		v_{t}=\left(-\frac{1}{a}\sqrt{\frac{b^{2}}{4}-ac}\right)\frac{1-k_{0}\exp\left(-2t\sqrt{\frac{b^{2}}{4}-ac}\right)}{1+k_{0}\exp\left(-2t\sqrt{\frac{b^{2}}{4}-ac}\right)}-\frac{b}{2a}~\text{ for }~t\in[0,T],
	\end{align}
	with $a:=-\beta^{2}\big(1+\rho_{rS}^{2}\big(1-\rho_{rS}^{2}\big)^{-1}\big),~b:=-2\kappa_{R}-2\beta\sigma_{R}\rho_{RS}+2\sigma_{R}\rho_{Rr}\beta\rho_{rS}\big(1-\rho_{rS}^{2}\big)^{-1/2}$, $c:=\sigma_{R}^{2}\big(1-\rho_{RS}^{2}-\rho_{Rr}^{2}\big)$ and $k_{0}:=\big(1+\frac{b}{2}\big(\frac{b^{2}}{4}-ac\big)^{-1/2}\big)\big(1-\frac{b}{2}\big(\frac{b^{2}}{4}-ac\big)^{-1/2}\big)^{-1}$.
	
	Moreover,
	\begin{align}\label{Boundedness of the conditional variance}
		0\le v_{t}<-\frac{1}{a}\sqrt{\frac{b^{2}}{4}-ac}-\frac{b}{2a} ~(=:v_{\infty})~\text{ for all }~t\in[0,T].
	\end{align}
\end{lemm}
\begin{proof}
	To check that $v$ given by \eqref{Solution of the variance equation_Expression} solves \eqref{Variance equation}, it suffices to differentiate $v$ and to compare the obtained expression with the right side of \eqref{Variance equation} for $v$ as in \eqref{Solution of the variance equation_Expression}. Uniqueness follows from the uniqueness of a solution to a Riccati equation. Observe that $a<0,c>0$. Then, $\big|\frac{b}{2}\big|<\sqrt{\frac{b^{2}}{4}-ac}$ and $k_{0}>0$. Using the expression of $v$, we obtain
	\begin{align}\label{Sign of the derivative of the condition variance}
		\frac{\mathrm{d}v_{t}}{\mathrm{d}t}&=-\frac{4}{a}\left(\frac{b^{2}}{4}-ac\right)k_{0}\frac{\exp\left(-2t\sqrt{\frac{b^{2}}{4}-ac}\right)}{\left(1+k_{0}\exp\left(-2t\sqrt{\frac{b^{2}}{4}-ac}\right)\right)^{2}}>0~\text{ for }~t\in[0,T].
	\end{align}
	In addition, $\lim\limits_{t\to\infty}v_{t}=-\frac{1}{a}\sqrt{\frac{b^{2}}{4}-ac}-\frac{b}{2a}=:v_{\infty}$. Hence, $v_{0}\le v_{t}<v_{\infty}$ for all $t\in[0,T]$.
\end{proof}

From the expression of $\eta$ just below \eqref{Filtered wealth}, we have
\begin{align}\label{Expression of eta^2}
	\|\eta_{t}\|^{2}=\frac{\beta^{2}}{1-\rho_{rS}^{2}}m_{t}^{2}-2\frac{\beta\rho_{rS}\phi_{B}}{\sqrt{1-\rho_{rS}^{2}}}m_{t}+\frac{\phi_{B}^{2}}{1-\rho_{rS}^{2}}~\text{ for }~t\in[0,T].
\end{align}

Hence, using the fact that $(m_{t})_{t\in[0,T]}$ is an OU process (see Proposition \ref{Conditional mean-variance pair}), Equation~\eqref{Boundedness of the conditional variance} and the convex inequality
\begin{align}\label{Convex inequality}
	\big(\sum_{i=1}^{\ell}a_{i}\big)^{q}\le\ell^{q-1}\sum_{i=1}^{\ell}a_{i}^{q}~\text{ for }~q\ge1,a_{i}>0,~i\in\{1,\cdots,\ell\}, 
\end{align}
we deduce that
\begin{align}\label{Integrability of eta and m_power}
	\mathbb{E}\Big[\int_{0}^{T}\|\eta_{s}\|^{2q}\mathrm{d}s\Big]+\mathbb{E}\Big[\int_{0}^{T}m_{s}^{q}\mathrm{d}s\Big]+\mathbb{E}\Big[\int_{0}^{T}r_{s}^{q}\mathrm{d}s\Big]+\mathbb{E}\Big[\exp\Big(-q\int_{0}^{T}m_{s}\mathrm{d}s\Big)\Big]<\infty,~q\ge0.
\end{align}

Now, using Cauchy-Schwarz inequality, we have 
\begin{align}
	\mathbb{E}\big[H_{T}^{q}\big]&\le\Big(\mathbb{E}\Big[\exp\Big((2q^{2}-q)\int_{0}^{T}\|\eta_{s}\|^{2}\mathrm{d}s\Big)\Big]\Big)^{2}\Big(\mathbb{E}\Big[\mathcal{E}\Big(\int-2q\eta^{\top}\mathrm{d}I\Big)_{s}\Big]\Big)^{2},\label{Integrability of HT_power}\\
	\mathbb{E}\Big[\int_{0}^{T}\varphi_{s}^{q}\mathrm{d}s\Big]&\le\Big(\mathbb{E}\Big[\exp\Big(\frac{2q^{2}+(\gamma-1)q}{\gamma^{2}}\int_{0}^{T}\|\eta_{t}\|^{2}\mathrm{d}t\Big)\Big]\Big)^{\frac{1}{2}}\notag\\
	&\phantom{Xx}\times\Big(\mathbb{E}\Big[\int_{0}^{T}\mathcal{E}\Big(\frac{2q}{\gamma}\int_{}^{}\eta^{\top}\mathrm{d}I\Big)_{t}\mathrm{d}t\Big]\Big)^{\frac{1}{2}}\widetilde{x}^{q'}e^{2q'\big(-\frac{\delta^{\psi}\psi}{\psi-1}+\frac{\delta\theta}{1-\gamma}\big)T}T^{\frac{1}{2}}.\label{Integrability of varphi_power}
\end{align}
and, using \eqref{Convex inequality}, we obtain
\begin{align}\label{Integrability of alpha_power}
	\mathbb{E}\Big[\int_{0}^{T}|\alpha_{s}|^{q}\mathrm{d}s\Big]&\le3^{q-1}\Big(T\Big|\frac{\delta^{\psi}}{\psi-1}-\frac{\delta\theta}{1-\gamma}\Big|^{q}+\mathbb{E}\Big[\int_{0}^{T}|r_{s}|^{q}\mathrm{d}s\Big]+\frac{1}{2^{q}\gamma^{q}}\mathbb{E}\Big[\int_{0}^{T}\|\eta_{s}\|^{2q}\mathrm{d}s\Big]\Big)<\infty.
\end{align}

Using \cite[Sect.~3.2.2.1 on p.64]{alos2008malliavin} and the product rule, we have
\begin{align}
	D_{t}\big(m_{s}\big)&=\mathds{1}_{[0,s]}(t)e^{\kappa_{R}(t-s)} \begin{pmatrix}
		\sigma_{R}\rho_{RS}+\beta v_{t}\\ \sigma_{R}\rho_{Rr}-\rho_{rS}\beta(1-\rho_{rS}^{2})^{-\frac{1}{2}}v_{t}
	\end{pmatrix}=:\mathds{1}_{[0,s]}(t)\begin{pmatrix}
		D_{t}^{(1)}\big(m_{s}\big)\\ D_{t}^{(2)}\big(m_{s}\big)
	\end{pmatrix}\label{E_B1}\\
	D_{t}\big(\eta_{s}\big)
	&=\mathds{1}_{[0,s]}(t)\begin{pmatrix}
		\beta D_{t}^{(1)}\big(m_{s}\big) & \beta D_{t}^{(2)}\big(m_{s}\big)\\ -\beta\rho_{rS}\big(1-\rho_{rS}^{2}\big)^{-\frac{1}{2}}D_{t}^{(1)}\big(m_{s}\big) & -\beta\rho_{rS}\big(1-\rho_{rS}^{2}\big)^{-\frac{1}{2}}D_{t}^{(2)}\big(m_{s}\big)
	\end{pmatrix}\label{E_B2}\\
	D_{t}\big(r_{s}\big)&=\mathds{1}_{[0,s]}(t)\sigma_{r}e^{\kappa_{r}(t-s)}\begin{pmatrix}
		\rho_{rS}\\ \sqrt{1-\rho_{rS}^{2}}
	\end{pmatrix}=:\mathds{1}_{[0,s]}(t)\begin{pmatrix}
		D_{t}^{(1)}\big(r_{s}\big)\\ D_{t}^{(2)}\big(r_{s}\big)
	\end{pmatrix}\label{E_B3}\\
	D_{t}\big(e^{-\int_{0}^{T}r_{s}\mathrm{d}s}\big)&=-e^{-\int_{0}^{T}r_{s}\mathrm{d}s}\int_{t}^{T}D_{t}(r_{s})\mathrm{d}s\label{E_B4}\\
	D_{t}\big(\alpha_{s}\big)
	&=\Big(\frac{\delta^{\psi}}{\psi-1}+r_{s}+\frac{1}{2\gamma}\|\eta_{s}\|^{2}-\frac{\delta\theta}{1-\gamma}\Big)D_{t}\Big(e^{-\int_{0}^{T}r_{s}\mathrm{d}s}\Big)\notag\\
	&\phantom{x}+e^{-\int_{0}^{s}r_{u}\mathrm{d}u}\Big(D_{t}\big(r_{s}\big)+\frac{\beta^{2}}{\gamma(1-\rho_{rS}^{2})}m_{s}D_{t}\big(m_{s}\big)-\frac{\beta\rho_{rS}\phi_{B}}{\gamma\sqrt{1-\rho_{rS}^{2}}}D_{t}\big(m_{s}\big)\Big).\label{E_B5}
\end{align}
Using Young inequality, \eqref{Boundedness of the conditional variance}, \eqref{Integrability of eta and m_power}, \eqref{E_B1}, \eqref{E_B3}, \eqref{E_B4} and \eqref{E_B5} we obtain
\begin{align}\label{Integrability of Malliavin derivatives_r and alpha_power}
	\mathbb{E}\Big[\int_{0}^{T}\big\|D_{t}\big(e^{-\int_{0}^{T}r_{s} \mathrm{d}s}\big)\big\|^{q}\mathrm{d}t\Big]+\mathbb{E}\Big[\int_{0}^{T}\Big(\int_{t}^{T}\|D_{t}\big(\alpha_{s}\big)\|^{q}\mathrm{d}s\Big)\mathrm{d}t\Big]<\infty\text{ for all }q>1.
\end{align}
Besides, using the expression of $\eta$ (see right below \eqref{Filtered wealth}) and the product rule, we have
\begin{align*}
	D_{t}\big(\varphi_{s}\big)
	&=\varphi_{s}\left(\frac{\gamma-1}{2\gamma^{2}}\int_{t}^{s}\begin{pmatrix}
		\big(2\frac{\beta^{2}}{1-\rho_{rS}^{2}}m_{u}-2\frac{\beta\rho_{rS}\phi_{B}}{\sqrt{1-\rho_{rS}^{2}}}\big)D_{t}^{(1)}\big(m_{u}\big)\\ \big(2\frac{\beta^{2}}{1-\rho_{rS}^{2}}m_{u}-2\frac{\beta\rho_{rS}\phi_{B}}{\sqrt{1-\rho_{rS}^{2}}}\big)D_{t}^{(2)}\big(m_{u}\big)
	\end{pmatrix}\mathrm{d}u\right.\\
	&\phantom{XXxxx}\left.+\frac{1}{\gamma}\eta_{t}+\int_{t}^{s}\begin{pmatrix}
		\frac{\beta}{\gamma}D_{t}^{(1)}\big(m_{u}\big) & 0\\ 0 & -\frac{\beta\rho_{rS}}{\gamma\sqrt{1-\rho_{rS}^{2}}}D_{t}^{(2)}\big(m_{u}\big)
	\end{pmatrix}\mathrm{d}I_{u}\right)\\
	D_{t}\left(H_{T}\right)
	&=H_{T}\left(-\frac{1}{2}\int_{t}^{T}\begin{pmatrix}
		\big(2\frac{\beta^{2}}{1-\rho_{rS}^{2}}m_{s}-2\frac{\beta\rho_{rS}\phi_{B}}{\sqrt{1-\rho_{rS}^{2}}}\big)D_{t}^{(1)}\big(m_{s}\big)\\ \big(2\frac{\beta^{2}}{1-\rho_{rS}^{2}}m_{s}-2\frac{\beta\rho_{rS}\phi_{B}}{\sqrt{1-\rho_{rS}^{2}}}\big)D_{t}^{(2)}\big(m_{s}\big)
	\end{pmatrix}\mathrm{d}s\right.\\
	&\phantom{XXxxx}\left.-\eta_{t}-\int_{t}^{T}\begin{pmatrix}
		\frac{\beta}{\gamma}D_{t}^{(1)}\big(m_{s}\big) & 0\\ 0 & -\frac{\beta\rho_{rS}}{\gamma\sqrt{1-\rho_{rS}^{2}}}D_{t}^{(2)}\big(m_{s}\big)
	\end{pmatrix}\mathrm{d}I_{s}\right).
\end{align*}

Using successively Young inequality, Minkowski inequality, Jensen inequality, \\Burkholder–Davis–Gundy (BDG) inequality and the convex inequality \eqref{Convex inequality}, we have
\begin{align}\label{Integrability of Malliavin derivatives_H_T}
	&\mathbb{E}\Big[\int_{0}^{T}\|D_{t}\big(H_{T}\big)\|^{q}\mathrm{d}t\Big]\notag\\
	&\le\mathbb{E}\Big[\int_{0}^{T}\Big(\Big|\frac{\beta^{2}}{1-\rho_{rS}^{2}}m_{s}\Big|^{2q(q+1)}+\big|D_{t}^{(1)}\big(m_{s}\big)\big|^{2q(q+1)}+\big|D_{t}^{(2)}\big(m_{s}\big)\big|^{2q(q+1)}\Big)\mathrm{d}t\Big]\notag\\
	&\phantom{X}+T\mathbb{E}\big[H_{T}^{q+1}\big]+T\Big|\frac{\beta\rho_{rS}\phi_{B}}{\sqrt{1-\rho_{rS}^{2}}}\Big|^{2q(q+1)}.
\end{align}
Similarly,
\begin{align}\label{Integrability of Malliavin derivatives_varphi_power}
	&\mathbb{E}\Big[\int_{0}^{T}\Big(\int_{t}^{T}\|D_{t}\big(\varphi_{s}\big)\|^{q}\mathrm{d}s\Big)\mathrm{d}t\Big]\notag\\
	&\le T\mathbb{E}\Big[\int_{0}^{T}\Big(\Big|\frac{\beta^{2}}{1-\rho_{rS}^{2}}m_{s}\Big|^{2q(q+1)}+\big|D_{t}^{(1)}\big(m_{s}\big)\big|^{2q(q+1)}+\big|D_{t}^{(2)}\big(m_{s}\big)\big|^{2q(q+1)}\Big)\mathrm{d}t\Big]\notag\\
	&\phantom{X}+T\mathbb{E}\Big[\int_{0}^{T}\varphi_{s}^{q+1}\mathrm{d}s\Big]+T^{2}\Big|\frac{\beta\rho_{rS}\phi_{B}}{\sqrt{1-\rho_{rS}^{2}}}\Big|^{2q(q+1)}
\end{align}
Hence, using \eqref{Boundedness of the conditional variance}, \eqref{Integrability of eta and m_power}, \eqref{E_B1}, \eqref{E_B3}, \eqref{E_B4} and \eqref{E_B5}, we deduce that
\begin{align}
	&\mathbb{E}\big[H_{T}^{q+1}\big]<\infty\Longrightarrow\mathbb{E}\Big[\int_{0}^{T}\|D_{t}\big(H_{T}\big)\|^{q}\mathrm{d}t\Big]<\infty\label{Finiteness of the Malliavin derivative of H_T}\\
	&\text{and }~\mathbb{E}\Big[\int_{0}^{T}\varphi_{s}^{q+1}\mathrm{d}s\Big]<\infty\Longrightarrow\mathbb{E}\Big[\int_{0}^{T}\Big(\int_{t}^{T}\|D_{t}\big(\varphi_{s}\big)\|^{q}\mathrm{d}s\Big)\mathrm{d}t\Big]<\infty.\label{Finiteness of the Malliavin derivative of varphi}
\end{align}

For $(i)$, $K_{T}e^{-\int_{0}^{T}r_{s}\mathrm{d}s}+\int_{0}^{T}\alpha_{s}\varphi_{s}\mathrm{d}s\in\mathbb{D}^{1,2}$, $H_{T}\Big(K_{T}e^{-\int_{0}^{T}r_{s}\mathrm{d}s}+\int_{0}^{T}\alpha_{s} \varphi_{s}\mathrm{d}s\Big)\in\mathbb{D}^{1,2}$.
\begin{align}\label{E1_Proposition 4.1}
	&\mathbb{E}\Big[\Big|K_{T}e^{-\int_{0}^{T}r_{s}\mathrm{d}s}+\int_{0}^{T}\alpha_{s}\varphi_{s}\mathrm{d}s\Big|^{2}\Big]\notag\\
	&\le2\mathbb{E}\Big[\Big|K_{T}e^{-\int_{0}^{T}r_{s}\mathrm{d}s}\Big|^{2}\Big]+2\mathbb{E}\Big[H_{T}^{-\frac{1}{2}}H_{T}^{\frac{1}{2}}\Big(\int_{0}^{T}\alpha_{s}\varphi_{s}\mathrm{d}s\Big)^{2}\Big]\notag\\
	&\le2\mathbb{E}\Big[\Big|K_{T}e^{-\int_{0}^{T}r_{s}\mathrm{d}s}\Big|^{2}\Big]+2\mathbb{E}\big[H_{T}^{-1}\big]+2T\mathbb{E}^{\mathbb{Q}_{0}}\Big[\Big(\int_{0}^{T}|\alpha_{s}\varphi_{s}|^{2}\mathrm{d}s\Big)^{2}\Big]\notag\\
	&<\infty,
\end{align}
where the first inequality comes from the convex inequality \eqref{Convex inequality}, the second inequality follows from Young inequality, Jensen inequality, and the last inequality holds due to Appendix \ref{Appendix A}, and Assumptions~\ref{Change of measure}$.(i)$ and \ref{Conditions to be Malliavin differentiable}$.(i)$.
\begin{align}\label{E2_Proposition 4.1}
	&\mathbb{E}\Big[\int_{0}^{T}\big\|D_{t}\Big(K_{T}e^{-\int_{0}^{T}r_{s}\mathrm{d}s}\Big)+D_{t}\Big(\int_{0}^{T}\alpha_{s}\varphi_{s}\mathrm{d}s\Big)\big\|^{2}\mathrm{d}t\Big]\notag\\
	&\le\mathbb{E}\Big[\int_{0}^{T}\big\|D_{t}\Big(K_{T}e^{-\int_{0}^{T}r_{s}\mathrm{d}s}\Big)\big\|^{2}\mathrm{d}t\Big]+T\mathbb{E}\Big[\int_{0}^{T}\Big(\int_{t}^{T}\big\|\alpha_{s}D_{t}\big(\varphi_{s}\big)+\varphi_{s}D_{t}\big(\alpha_{s}\big)\big\|^{2}\mathrm{d}s\Big)\mathrm{d}t\Big]\notag\\
	&\le\mathbb{E}\Big[\int_{0}^{T}\big\|D_{t}\Big(K_{T}e^{-\int_{0}^{T}r_{s}\mathrm{d}s}\Big)\big\|^{2}\mathrm{d}t\Big]+T^{2}\mathbb{E}\Big[\int_{0}^{T}\big(2\varphi_{s}^{3}+\alpha_{s}^{6}\varphi_{s}^{-6}\|D_{t}\big(\varphi_{s}\big)\|^{6}+\|D_{t}\big(\alpha_{s}\big)\|^{6}\big)\mathrm{d}s\Big]\notag\\
	&<\infty,
\end{align}
where the first inequality holds due to Minkowski inequality, Jensen inequality, the second inequality comes from the triangle inequality and Young inequality, and the last inequality follows from Equations \eqref{Integrability of varphi_power} (for $q=3$), \eqref{Integrability of eta and m_power}, \eqref{Integrability of alpha_power} and \eqref{Integrability of Malliavin derivatives_r and alpha_power}, and Assumptions~\ref{Change of measure}$.(i)$ and \ref{Conditions to be Malliavin differentiable}$.(i)$.
\begin{align}\label{E4_Proposition 4.1}
	&\mathbb{E}\Big[H_{T}^{2}\Big|K_{T}e^{-\int_{0}^{T}r_{s}\mathrm{d}s}+\int_{0}^{T}\alpha_{s}\varphi_{s}\mathrm{d}s\Big|^{2}\Big]\notag\\
	&\le2\mathbb{E}\Big[H_{T}^{2}\Big|K_{T}e^{-\int_{0}^{T}r_{s}\mathrm{d}s}\Big|^{2}\Big]+2\mathbb{E}\Big[H_{T}^{2}\Big(\int_{0}^{T}\alpha_{s}\varphi_{s}\mathrm{d}s\Big)^{2}\Big]\notag\\
	&\le2\mathbb{E}\Big[H_{T}^{2}\Big|K_{T}e^{-\int_{0}^{T}r_{s}\mathrm{d}s}\Big|^{2}\Big]+2T\mathbb{E}\Big[H_{T}^{\frac{3}{2}}\Big(H_{T}^{\frac{1}{2}}\int_{0}^{T}|\alpha_{s}\varphi_{s}|^{2}\mathrm{d}s\Big)\Big]\notag\\
	&\le2\mathbb{E}\Big[H_{T}^{2}\Big|K_{T}e^{-\int_{0}^{T}r_{s}\mathrm{d}s}\Big|^{2}\Big]+2T\Big(\mathbb{E}\big[H_{T}^{3}\big] +\mathbb{E}\Big[H_{T}\Big(\int_{0}^{T}|\alpha_{s}\varphi_{s}|^{2}\mathrm{d}s\Big)^{2}\Big]\Big)\notag\\
	&=2\mathbb{E}\Big[H_{T}^{2}\Big|K_{T}e^{-\int_{0}^{T}r_{s}\mathrm{d}s}\Big|^{2}\Big]+2T\Big(\mathbb{E}\big[H_{T}^{3}\big] +\mathbb{E}^{\mathbb{Q}_{0}}\Big[\Big(\int_{0}^{T}|\alpha_{s}\varphi_{s}|^{2}\mathrm{d}s\Big)^{2}\Big]\Big)\notag\\
	&<\infty,
\end{align}
where the first inequality holds due to the convex inequality \eqref{Convex inequality}, the second inequality comes from Jensen inequality, the third follows from Young inequality, and the last inequality comes from Equations \eqref{Integrability of HT_power} (for $q=3$), Appendix \ref{Appendix A}, and Assumption \ref{Change of measure}$.(i)$ and \ref{Conditions to be Malliavin differentiable}$.(i)$.
\begin{align}\label{E5_Proposition 4.1}
	&\mathbb{E}\Big[\int_{0}^{T}\big\|D_{t}\Big(H_{T}K_{T}e^{-\int_{0}^{T}r_{s}\mathrm{d}s}\Big)+D_{t}\Big(H_{T}\int_{0}^{T}\alpha_{s} \varphi_{s}\mathrm{d}s\Big)\big\|^{2}\mathrm{d}t\Big]\notag\\
	&\le\mathbb{E}\Big[\int_{0}^{T}\big\|D_{t}\Big(H_{T}K_{T}e^{-\int_{0}^{T}r_{s}\mathrm{d}s}\Big)\big\|^{2}\mathrm{d}t\Big]+\mathbb{E}\Big[\int_{0}^{T}\Big\|\Big(\int_{0}^{T}\alpha_{s}\varphi_{s}\mathrm{d}s\Big)D_{t}\big(H_{T}\big)\Big\|^{2}\mathrm{d}t\Big]\notag\\
	&\phantom{X}+\mathbb{E}\Big[\int_{0}^{T}\Big\|H_{T}D_{t}\Big(\int_{0}^{T}\alpha_{s}\varphi_{s}\mathrm{d}s\Big)\Big\|^{2}\mathrm{d}t\Big]\notag\\
	&\le\mathbb{E}\Big[\int_{0}^{T}\big\|D_{t}\Big(H_{T}K_{T}e^{-\int_{0}^{T}r_{s}\mathrm{d}s}\Big)\big\|^{2}\mathrm{d}t\Big]+\mathbb{E}\Big[\int_{0}^{T}H_{T}^{\frac{1}{2}}\Big(\int_{0}^{T}\alpha_{s}\varphi_{s}\mathrm{d}s\Big)^{2}\Big\|H_{T}^{-\frac{1}{4}}D_{t}\big(H_{T}\big)\Big\|^{2}\mathrm{d}t\Big]\notag\\
	&\phantom{X}+T\mathbb{E}\big[H_{T}^{4}\big]+\mathbb{E}\Big[\int_{0}^{T}\Big\|D_{t}\Big(\int_{0}^{T}\alpha_{s}\varphi_{s}\mathrm{d}s\Big)\Big\|^{4}\mathrm{d}t\Big]\notag\\
	&\le\mathbb{E}\Big[\int_{0}^{T}\big\|D_{t}\Big(H_{T}K_{T}e^{-\int_{0}^{T}r_{s}\mathrm{d}s}\Big)\big\|^{2}\mathrm{d}t\Big]+T\mathbb{E}\Big[H_{T}\Big(\int_{0}^{T}\alpha_{s}\varphi_{s}\mathrm{d}s\Big)^{4}\Big]+2T\mathbb{E}\big[H_{T}^{4}\big]\notag\\
	&\phantom{X}+\mathbb{E}\Big[\int_{0}^{T}\Big\|H_{T}^{-1}D_{t}\big(H_{T}\big)\Big\|^{8}\mathrm{d}t\Big]+T^{4}\mathbb{E}\Big[\int_{0}^{T}\big(2\varphi_{s}^{5}+\alpha_{s}^{20}\varphi_{s}^{-20}\|D_{t}\big(\varphi_{s}\big)\|^{20}+\|D_{t}\big(\alpha_{s}\big)\|^{20}\big)\mathrm{d}s\Big]\notag\\
	&<\infty,
\end{align}
where the first inequality holds due to the triangle inequality and the product rule, the second inequality comes from Young inequality, the third inequality follows from the triangle inequality and Young inequality, and the last inequality follows from Equations \eqref{Integrability of HT_power} (for $q=4$), \eqref{Integrability of varphi_power} (for $q=5$), \eqref{Integrability of eta and m_power}, \eqref{Integrability of alpha_power} and \eqref{Integrability of Malliavin derivatives_r and alpha_power}, and Assumptions~\ref{Change of measure}$.(i)$ and \ref{Conditions to be Malliavin differentiable}$.(i)$.

For $(ii)$, we need to show that $\mathbb{E}\Big[\|\eta_{t}\|^{2}+\int_{0}^{T}\|D_{s}(\eta_{t})\|^{2}\mathrm{d}s\Big]<\infty$. This follows directly from Equations \eqref{Boundedness of the conditional variance}, \eqref{Integrability of eta and m_power}, \eqref{E_B1} and \eqref{E_B2}.

For $(iii)$-$(iv)$,
\begin{align}
	&\mathbb{E}^{\mathbb{Q}_{0}}\Big[\big|K_{T}e^{-\int_{0}^{T}r_{s}\mathrm{d}s}+\int_{0}^{T}\alpha_{s}\varphi_{s}\mathrm{d}s\big|\Big]\notag\\
	&+\mathbb{E}^{\mathbb{Q}_{0}}\Big[\int_{0}^{T}\Big(\big\|D_{t}\big(K_{T}e^{-\int_{0}^{T}r_{s}\mathrm{d}s}\big)\big\|^{2}+\big\|D_{t}\big(\int_{0}^{T}\alpha_{s}\varphi_{s}\mathrm{d}s\big)\big\|^{2}\Big)\mathrm{d}t\Big]\notag\\
	&\le\mathbb{E}^{\mathbb{Q}_{0}}\Big[\big|K_{T}e^{-\int_{0}^{T}r_{s}\mathrm{d}s}\big|+\int_{0}^{T}\big\|D_{t}\big(K_{T}e^{-\int_{0}^{T}r_{s}\mathrm{d}s}\big)\big\|^{2}\mathrm{d}t\Big]+\mathbb{E}^{\mathbb{Q}_{0}}\Big[\big|\int_{0}^{T}\alpha_{s}\varphi_{s}\mathrm{d}s\big|\Big]\notag\\
	&\phantom{X}+\mathbb{E}^{\mathbb{Q}_{0}}\Big[\int_{0}^{T}\big\|D_{t}\big(\int_{0}^{T}\alpha_{s}\varphi_{s}\mathrm{d}s\big)\big\|^{2}\mathrm{d}t\Big]\notag\\
	&\le\mathbb{E}^{\mathbb{Q}_{0}}\Big[\big|K_{T}e^{-\int_{0}^{T}r_{s}\mathrm{d}s}\big|+\int_{0}^{T}\big\|D_{t}\big(K_{T}e^{-\int_{0}^{T}r_{s}\mathrm{d}s}\big)\big\|^{2}\mathrm{d}t\Big]+T^{\frac{1}{2}}\Big(\mathbb{E}^{\mathbb{Q}_{0}}\Big[\Big(\int_{0}^{T}|\alpha_{s}\varphi_{s}|^{2}\mathrm{d}s\Big)^{2}\Big]\Big)^{\frac{1}{4}}\notag\\
	&\phantom{X}+\mathbb{E}\big[H_{T}^{2}\big]+T\mathbb{E}\Big[\int_{0}^{T}\big\|D_{t}\big(\int_{0}^{T}\alpha_{s}\varphi_{s}\mathrm{d}s\big)\big\|^{4}\mathrm{d}t\Big]\notag\\
	&<\infty,
\end{align}
where the first inequality holds by sub-additivity of the absolute value function, the second inequality comes from Young inequality, H{\"o}lder inequality and Jensen inequality, and the last inequality follows from Equation \eqref{E5_Proposition 4.1} (see its second inequality), Appendix \ref{Appendix A}, and Assumptions~\ref{Change of measure}$.(i)$ and \ref{Conditions to be Malliavin differentiable}$.(i)$.

To prove $(v)$, it suffices to show that
\begin{align}\label{E_(v)_1}
	\mathbb{E}^{\mathbb{Q}_{0}}\Big[\int_{0}^{T}\|D_{t}(\alpha_{t}\varphi_{t})\|^{2}\mathrm{d}t\Big]<\infty~\text{ and }~\mathbb{E}^{\mathbb{Q}_{0}}\Big[\int_{0}^{T}\|Z_{t}^{\top}D_{t}(\eta_{t})\|^{2}\mathrm{d}t\Big]<\infty.
\end{align}
Because the proof of the first inequality in \eqref{E_(v)_1} is on similar lines with the proof of \eqref{E2_Proposition 4.1}, it is omitted for brevity. It remains to show the second inequality in \eqref{E_(v)_1}. First, note that $\max\big(D_{t}^{(1)}\big(m_{s}\big),~D_{t}^{(2)}\big(m_{s}\big)\big)<\infty$ for all $t,s\in[0,T]$; thanks to Equations \eqref{E_B1} and
\eqref{Boundedness of the conditional variance}. Hence, using \eqref{E_B1}, \eqref{E_B2} and the fact that $\mathbb{E}^{\mathbb{Q}_{0}}\Big[\Big(\int_{0}^{T}\|Z_{t}\|^{2}\mathrm{d}t\Big)^{2}\Big]<\infty$ (see the proof of Proposition \ref{Existence result for the FBSDE}; there in $Z=\tilde{Z}$ and $q=2$), we have
\begin{align}
	\mathbb{E}^{\mathbb{Q}_{0}}\Big[\int_{0}^{T}\|Z_{t}^{\top}D_{t}(\eta_{t})\|^{2}\mathrm{d}t\Big]&\le\max\big(D_{t}^{(1)}\big(m_{s}\big),~D_{t}^{(2)}\big(m_{s}\big)\big)\mathbb{E}^{\mathbb{Q}_{0}}\Big[\int_{0}^{T}\|Z_{t}\|^{2}\mathrm{d}t\Big]\notag\\
	&\le\max\big(D_{t}^{(1)}\big(m_{s}\big),~D_{t}^{(2)}\big(m_{s}\big)\big)\Big(\mathbb{E}^{\mathbb{Q}_{0}}\Big[\Big(\int_{0}^{T}\|Z_{t}\|^{2}\mathrm{d}t\Big)^{2}\Big]\Big)^{\frac{1}{2}}\notag\\
	&<\infty.
\end{align}
This concludes the proof of Proposition \ref{Auxilliary results_Malliavin differentiability}.

\section{Proof of Proposition \ref{Sufficient conditions for the assumptions to hold}}\label{Appendix D}
Recall that expectations are taken under the observation filtration $\mathbb{F}^{r,S}$. We first state and prove two intermediate results (Lemmas \ref{Exponential moment of the square of an OU process} and \ref{Comparison theorem for Ricatti equations}) on which the proof of Proposition \ref{Sufficient conditions for the assumptions to hold} will rely on. Lemma \ref{Exponential moment of the square of an OU process} gives sufficient conditions for the non-explosion of the exponential moments of the square of an OU process with time-dependent coefficients and Lemma \ref{Comparison theorem for Ricatti equations} gives a comparison result for some Riccati equations.
\begin{lemm}\label{Exponential moment of the square of an OU process}
	For $\zeta=2q_{max}\beta^{2}(1-\rho_{rS}^{2})^{-1}>0$, let $\Delta_{max}:=2b_{max}\zeta-\kappa_{R}^{2}$. If $\Delta_{max}\le0$ or $\Delta_{max}>0,~T<\left(\textsl{pi}-\arctan(\sqrt{\Delta_{max}}/\kappa_{R})\right)/\sqrt{\Delta_{max}}$ hold, then $\mathbb{E}\left[\exp\left(\zeta\int_{0}^{T}m_{t}^{2}\mathrm{d}t\right)\right]<\infty$.
\end{lemm}
\begin{proof}
	Define $u(t,m):=\mathbb{E}\left[\exp\left(\zeta\int_{t}^{T}m_{s}^{2}\mathrm{d}s\right)\big|m_{t}=m\right]$. Then $u$ satisfies the backward Feynman–K{\v a}c partial differential equation (PDE):
	\begin{align}
		\frac{\partial u}{\partial t}-\kappa_{R}m\frac{\partial u}{\partial m} +\frac{1}{2}\sigma_{m}^{2}(t)\frac{\partial^2 u}{\partial m^2}+\zeta m^{2}u=0,~\text{ with }~u(T,m)=1.
	\end{align}
	We make the exponential–quadratic ansatz $u(t,m)=\exp\big(g(t)m^{2}+B(t)\big)$, with $g(t)=0,B(T)=0$. Hence, $u_{t}=\big(g'(t)m^{2}+B'(t)\big)u,u_{m}=2g(t)mu,~u_{mm}=(2g(t)+4g^{2}(t)m^{2})u$ and we have
	\begin{align}
		&\big(g'(t)-2\kappa_{R}g(t)+2\sigma_{m}^{2}(t)g^{2}(t)+\zeta\big)m^{2}+B'(t)+\sigma_{m}^{2}(t)g(t)=0~\text{ for all }~m\in\mathbb{R}.
	\end{align}
	Hence,
	\begin{align}
		g'(t)=-2\sigma_{m}^{2}(t)g^{2}(t)+2\kappa_{R}g(t)-\zeta~\text{ and }~ B'(t)=-\sigma_{m}^{2}(t)g(t).
	\end{align}
	Using Lemma \ref{Comparison theorem for Ricatti equations} we have $0\le\frac{\zeta}{2\kappa_{R}}\left(\exp\big(2\kappa_{R}(T-t)\big)-1\right)\le g(t)\le g_{2}(t)$ and $B(t)\le0$, with $\frac{\zeta}{2\kappa_{R}}\big(\exp\big(2\kappa_{R}(T-t)\big)-1\big)=g_{3}(t)$ for all $t\in[0,T]$.
	
	Therefore, from the exponential-quadratic ansatz we obtain
	\begin{align}\label{E1_Lemma B3}
		\mathbb{E}\left[\exp\left(\zeta\int_{0}^{T}m_{t}^{2}\mathrm{d}t\right)\right]\le\exp\big(g_{2}(0)m^{2}\big).
	\end{align}
	
	Now, we solve the Riccati equation satisfied by $g_{2}$. We consider the transformation $g_{2}(t)=\frac{1}{2b_{max}^{2}}\frac{g_{4}'(t)}{g_{4}(t)}$. Then $g_{2}'(t)=\frac{g_{4}''(t)g_{4}(t)-(g_{4}'(t))^2}{2b_{max}^{2}g_{4}^{2}(t)}$. Hence, $g_{4}$ satisfies the linear ODE $g_{4}''=2\kappa_{R}g_{4}'-2b_{max}^{2}\zeta g_{4}$. Thus,
	\begin{align}
		g_{4}(t)&=k_{1}e^{(\kappa_{R}+\sqrt{-\Delta_{max}})t}+k_{2}e^{(\kappa_{R}-\sqrt{-\Delta_{max}})t},~\text{with } \Delta_{max}=2b_{max}^{2}\zeta-\kappa_{R}^{2}.
	\end{align}
	Hence,
	\begin{align}\label{E2_Lemma B3}
		g_{2}(t)&=\frac{\kappa_{R}\big(k_{1}e^{(\kappa_{R}+\sqrt{-\Delta_{max}})t}+k_{2}e^{(\kappa_{R}-\sqrt{-\Delta_{max}})t}\big)}{2b_{max}^{2}\big(k_{1}e^{(\kappa_{R}+\sqrt{-\Delta_{max}})t}+k_{2}e^{(\kappa_{R}-\sqrt{-\Delta_{max}})t}\big)}\notag\\
		&+\frac{\sqrt{-\Delta_{max}}\big(k_{1}e^{(\kappa_{R}+\sqrt{-\Delta_{max}})t}-k_{2}e^{(\kappa_{R}-\sqrt{-\Delta_{max}})t}\big)}{2b_{max}^{2}\big(k_{1}e^{(\kappa_{R}+\sqrt{-\Delta_{max}})t}+k_{2}e^{(\kappa_{R}-\sqrt{-\Delta_{max}})t}\big)}.
	\end{align}
	Applying the boundary condition $g_{2}(T)=0$ to fix the constants $k_{1},k_{2}$ we obtain
	\begin{align}
		g_{2}(0)&=\frac{\zeta\sinh(T\sqrt{-\Delta_{max}})}{2\left(\sqrt{-\Delta_{max}}\cosh(T\sqrt{-\Delta_{max}})+\kappa_{R}\sinh(T\sqrt{-\Delta_{max}})\right)}.
	\end{align}
	
	Next, we discuss the finiteness of $g_{2}(0)$. We obtain the following situations.\\
	Case $1$: For $\Delta_{max}<0$, the denominator of the fraction on the right side of \eqref{E2_Lemma B3} does not vanish. Then $g_{2}(0)<\infty$.\\
	Case $2$: For $\Delta_{max}=0$, the denominator as well as the numerator of the fraction on the right side of \eqref{E2_Lemma B3} vanishes. However, $g_{2}(0)=\frac{1}{2}\zeta T\big(1+\kappa_{R}T\big)^{-1}<\infty$.\\
	Case $3$: For $\Delta_{max}>0$, the denominator of the fraction on the right side of \eqref{E2_Lemma B3} does not vanish for all $T$ smaller than a critical value $T_{c}$. Indeed, using the facts that $\sqrt{-\Delta_{max}}=i\sqrt{\Delta_{max}}$, $\sinh(iT\sqrt{\Delta_{max}})=i\sin(T\sqrt{\Delta_{max}})$ and $\cosh(iT\sqrt{\Delta_{max}})=\cos(T\sqrt{\Delta_{max}})$ we have
	\begin{align}\label{E3_Lemma B3}
		g_{2}(0)&=\frac{\zeta\sin(T\sqrt{\Delta_{max}})}{2\left(\sqrt{\Delta_{max}}\cos(T\sqrt{\Delta_{max}})+\kappa_{R}\sin(T\sqrt{\Delta_{max}})\right)}.
	\end{align}
	Finding the first positive $T$ such that $\sqrt{\Delta_{max}}\cos(T\sqrt{-\Delta_{max}})+\kappa_{R}\sin(T\sqrt{\Delta_{max}})=0$ is equivalent to find the smallest $T>0$ satisfying $\tan(T\sqrt{-\Delta_{max}})=\frac{\sqrt{\Delta_{max}}}{\kappa_{R}}$. If we denote by $T_{c}$ such value, then $T_{c}=\frac{1}{\sqrt{\Delta_{max}}}\left( \textsl{pi}-\arctan\left(\frac{\sqrt{\Delta_{max}}}{\kappa_{R}}\right) \right)$. Hence, $g_{2}(0)<\infty$ for all $T<T_{c}$.
	
	Finally, using \eqref{E1_Lemma B3}  and the results in Cases $1$-$3$ we conclude the proof.
\end{proof}
\begin{lemm}\label{Comparison theorem for Ricatti equations}
	For $v_{t},t\in[0,T]$, defined as in Proposition \ref{Conditional mean-variance pair}, let $\sigma_{m}^{2}(t):=\big(\sigma_{R}\rho_{RS}+\beta v_{t}\big)^{2}+\big(\sigma_{R}\rho_{Rr} -\rho_{rS}\beta(1-\rho_{rS}^{2})^{-\frac{1}{2}}v_{t}\big)^{2},~t\in[0,T]$, and $b_{max}^{2}:=\max_{t\in[0,T]}\sigma_{m}^{2}(t)$. If $g_{1},g_{2}$ and $g_{3}$ are solutions on $[0,T]$ of the ordinary equations
	\begin{align*}
		&g_{1}'(t)=-2\sigma_{m}^{2}(t)g_{1}^{2}(t)+2\kappa_{R}g_{1}(t)-\zeta,\quad g_{2}'(t)=-2b_{max}^{2}g_{2}^{2}(t)+2\kappa_{R}g_{2}(t)-\zeta\\
		&\text{and }~g_{3}'(t)=2\kappa_{R}g_{3}(t)-\zeta
	\end{align*}
	with $g_{1}(T)=g_{2}(T)=g_{3}(T)$, then $g_{3}(t)\le g_{1}(t)\le g_{2}(t)$ for all $t\in[0,T]$.
\end{lemm}
\begin{proof}
	The proof follows from theorem $4.1.4$ (on p.$185$) in \cite{abou2012matrix}.
\end{proof}
We can now confirm Proposition \ref{Sufficient conditions for the assumptions to hold}.
\begin{proof}[Proof of Proposition \ref{Sufficient conditions for the assumptions to hold}]
	Recall that $\zeta:=2q_{max}\beta^{2}(1-\rho_{rS}^{2})^{-1}$.
	\begin{align*}
		&\mathbb{E}\Big[\exp\Big(q_{max}\int_{0}^{T}\|\eta_{s}\|^{2}\mathrm{d}s\Big)\Big]\notag\\
		&=\mathbb{E}\Big[\exp\Big(q_{max}\int_{0}^{T}\Big(\frac{\beta^{2}}{1-\rho_{rS}^{2}}m_{s}^{2}-2\frac{\beta\rho_{rS}\phi_{B}}{\sqrt{1-\rho_{rS}^{2}}}m_{s}+\frac{\phi_{B}^{2}}{1-\rho_{rS}^{2}}\Big)\mathrm{d}s\Big)\Big]\notag\\
		&\le e^{q_{max}\phi_{B}^{2}(1-\rho_{rS}^{2})^{-1}T}\Big(\mathbb{E}\Big[e^{-4q_{max}\beta\rho_{rS}\phi_{B}(1-\rho_{rS}^{2})^{-1/2}\int_{0}^{T}m_{s}\mathrm{d}s}\Big]\Big)^{\frac{1}{2}}\Big(\mathbb{E}\Big[e^{2q_{max}\beta^{2}(1-\rho_{rS}^{2})^{-1}\int_{0}^{T}m_{s}^{2}\mathrm{d}s}\Big]\Big)^{\frac{1}{2}}\notag\\
		&<\infty,
	\end{align*}
	where the first inequality holds due to Cauchy-Schwarz inequality and the last inequality comes from \eqref{Integrability of eta and m_power}, \eqref{E1_Lemma B3}, Lemma \ref{Exponential moment of the square of an OU process} and the fact that $2q_{max}\beta^{2}(1-\rho_{rS}^{2})^{-1}=\zeta$.
\end{proof}

\section*{Acknowledgments}
I would like to acknowledge fruitful discussions with Olivier Menoukeu Pamen. I thank Thorsten Schmidt for his careful reading of the first version. I am grateful to an anonymous referee for constructive comments and suggestions that led to a substantial enhancement of the paper.

\bibliographystyle{siamplain}
\bibliography{references}

\end{document}